\documentclass{article}
\usepackage[top=30truemm,bottom=30truemm,left=25truemm,right=25truemm]{geometry}
\usepackage{amsmath,amsthm,amssymb,amsfonts,mathtools,accents,bm,mathrsfs,bibentry,lscape,multirow,natbib}
\usepackage{enumerate}
\usepackage{tabu}
\newtheorem{theorem}{Theorem}
\newtheorem{lemma}{Lemma}

\title{A semiparametric probability distribution estimator of sample maximums}
\author{
Taku MORIYAMA\\
School of Data Science, Yokohama City University}
\date{}

\begin{document}
\maketitle

\begin{abstract}
This study proposes a computationally efficient semiparametric distribution estimator, which is a slight modification of the naive mixture proposed by Schuster and Yakowitz (1985) and Olkin and Spiegelman (1987). The proposed method is applied to probability distribution estimation of a sample maximum. Two approaches for the sample maximum distribution estimation, one based on extreme value theory and the other on nonparametric smoothing, exist; however, theoretical and numerical properties of the two approaches are known to heavily depend on {the} case and greatly differ. This study demonstrates that the semiparametric mixture distribution estimator{s} ha{ve} good properties of both approaches. The cross-validation method is proposed for the mixing ratio selection {for the proposed mixture distribution estimator. T}he result of simulation experiments and {three case studies} are reported.
\end{abstract}
{\it Keywords:} Cross-validation; distribution estimator; sample maximum; semiparametric estimation

\section{Introduction of the mixture estimators}

It is known that there are cases in which nonparametric estimators of $F$ do not work, such as when the distribution is not smooth or too variational, when the rank of the underlying structure is too high, or when the interest is in the tail of the distribution (\citealt{smith1987estimating}). In these cases the mixture approach is possibly more effective than the naive nonparametric one. 

Recently, a semiparametric approach was developed. \cite{hjort1996locally} and \cite{loader1996local} proposed a locally parametric nonparametric density estimator obtained as the minimizer of the local kernel-smoothed likelihood. \cite{hjort1995nonparametric} proposed a parametrically guided nonparametric density estimator. The latter approach was extended to a case of censored data by \cite{talamakrouni2016parametrically}, and the asymptotic properties of the density and hazard function estimators were derived. 

Being intuitive, the semiparametrically mixing approach was proposed by \cite{schuster1985parametric} and \cite{olkin1987semiparametric}. The mixture distribution estimator has a target parametric distribution class and is given by
\begin{align*}
p F_{\widehat{\theta}}(x) + (1- p) \widehat{F}_{\widehat{h}}(x),
\end{align*}
where $F_{\widehat{\theta}}(x)$ is the parametric estimator with the MLE $\widehat{\theta}$ and $\widehat{F}_{\widehat{h}}(x)$ is the kernel distribution estimator, respectively. $0\le p \le 1$ is the mixing ratio parameter. If the underlying distribution $F$ belongs to the parametric class $\mathcal{F}_{\theta}:=\{F_{\theta} \vert \theta \in \Theta\}$, the mixture density estimator has $\sqrt{n}$ consistency under the model assumption (\citealt{olkin1987semiparametric}). \cite{faraway1990implementing} conducted a simulation study and investigated numerical properties. There exist the mixing ratio selectors based on pseudolikelihood approach (\citealt{olkin1987semiparametric}), minimized integrated squared error approach with a large parametric family of distributions (\citealt{rahman1997note}), bootstrap approach (\citealt{soleymani2014bootstrap}).

The naive kernel density estimator is seen as the locally fitting estimator to a uniform distribution (\citealt{hjort1996locally}). On the kernel cumulative distribution estimator it holds that
$$\widehat{F}_h(x) \to 
\begin{dcases}
\bar{F}(x) ~~~ & {\rm as} ~~~ h \to 0 \\
0.5 ~~~ & {\rm as} ~~~ h \to \infty,
\end{dcases}
$$
where $\bar{F}$ is the empirical distribution estimator. $\bar{F}(x)$ is unbiased and has the variance $(n^{-1} F(x) \{1- F(x)\})$. $0.5$ has zero variance and is biased unless $F$ is the uniform distribution on the whole real line. Thus, the bandwidth in the kernel distribution estimation controls the bias-variance trade-off. Therefore, in the mixture approach proposed by \cite{schuster1985parametric} and \cite{olkin1987semiparametric} the bandwidth also needs to be appropriately chosen; hence, the mixture density estimator consists of the two-steps. 

Thus, the two-steps algorithm is considered to aim at the uniform distribution and adjusts the mixing ratio in a sense (see also \cite{jones1993kernel} considering the large bandwidth). This study insists the first step bandwidth selection can be skipped for the sake of computational cost reduction without significant loss of accuracy and proposes a slightly modified estimator
\begin{align}
\widehat{F} (x; h) := q F_{\widehat{\theta}}(x) + (1- q) \widehat{F}_{\widehat{h}}(x),
\end{align}
where $0\le q \le 1$ is a designated function of $h$ satisfying
\begin{align}
q \to 
\begin{dcases}
0 ~~~ & {\rm as} ~~~ h \to 0 \\
1 ~~~ & {\rm as} ~~~ h \to \infty.
\end{dcases}
\end{align}
$h$ is also the smoothing parameter of the kernel cumulative distribution estimator. Therefore, the number of tuning parameter is one. Unlike the ordinary mixture approach, the proposed estimator does not have the secondary hyperparameter, but (1) satisfies
$$\widehat{F} (x; h) \to 
\begin{dcases}
\bar{F}(x) ~~~ & {\rm as} ~~~ h \to 0 \\
F_{\widehat{\theta}}(x) ~~~ & {\rm as} ~~~ h \to \infty.
\end{dcases}
$$
$h$ is expected to control the MSE trade-off directly. Thus, the proposed mixture estimator is expected to be computationally effective. 

Although one does not usually know the parametric class of the underlying distribution, there are a few cases in which we know the class of the approximating distribution. In the cases semiparametric approaches are possibly quite effective, and one of the cases is the sample maximum distribution (SMD) estimation. Examples of practical application include the evaluation of the sample maximum (\citealt{komukai1994requirements}; \citealt{kasai2016predicting}), outlier detection (\citealt{mittnik2001distribution}; \citealt{gbenro2020using}) and measuring future risk (\citealt{beirlant1996practical}; \citealt{resnick1997discussion}). 

\cite{moriyama2025application} demonstrated that the accuracy of the parametrically fitting estimator and the kernel type estimator are comparable and in some cases kernel type estimator does not work at all. Which should be used complexly depends on the case, and SMD estimation is considered to be an academically interesting case for suggesting the existence of an obvious limitation in the naive nonparametric estimator and the need for a new approach that is neither wholly nonparametric nor parametric. Motivated by the interest we investigate the mixture approach; however, due to the request for a large sample size in SMD estimation the computational load of the bandwidth and the mixing ratio selection is quite expensive. This study proposes an alternative approach, which determines the bandwidth and the mixing ratio simultaneously. The effectiveness of the proposed mixture approach is demonstrated. 

The rest of this paper is organized as follows. Basic properties of the proposed mixture estimator of SMD are given in Section 2. Simulation results on the naive mixture and the proposed mixture approaches are presented in Section 3. Section 4 reports {three} case studies, and the mixing parameter selection is discussed in Section 5. Section 6 concludes this study.

\section{Properties of the proposed mixture approach in SMD estimation}

Suppose $X_1, X_2, \cdots, X_n$ is the i.i.d. sequence, whose distribution function is $F$ and $F^{m}$ is the target distribution function. The nonparametric (kernel-type) estimator NE can be the kernel cumulative distribution the to the power of $m$ i.e.
$$(\widehat{F}_h)^m(x) := \left\{\frac{1}{n} \sum_{i=1}^n W\left(\frac{x- X_i}{h}\right) \right\}^m,$$
where $W$ is the distribution function of a symmetric density function $w$ satisfying $w(0) \neq 0$. 

The classical extreme value theory states that SMD $F^m$ is possibly approximated by the generalized extreme value distribution (GEV) $G_{\bm{\gamma}}$ for sufficiently large $m$. In that case SMD estimation means that the parametric model is asymptotically specified correctly. Let us assume $N:=n/m$ be an integer and $Y_1, \cdots, Y_N$ be the block maxima, where $Y_j:=\max\{X_{m(j-1)+1}, \cdots, X_{mj}\}$ for $j=1,\cdots,N$. Set $\widehat{\bm{\gamma}}$ is the maximum-likelihood estimator based on $\{Y_j\}_{j=1}^N$. 

The convergence rates of MSE of the parametrically fitting estimator PE $G_{\widehat{\bm{\gamma}}}(x)$ and NE $(\widehat{F}_h)^m(x)$ are given in Table {12} in the supplementary file, where 
\begin{align}
m \to \infty, ~\frac{m}{n} \to 0, ~ x \to x^*:=\sup({\rm supp}(f)) ~ {\rm and} ~ \ln F^m(x) =O(1)
\end{align}
as $n \to \infty$ (provided by \citealt{moriyama2025application}). The hyphen Table {12} means the assumption of the corresponding estimator is broken, and so the optimal convergence rate could not be specified. Since the approximation error of the GEV fitting becomes large as the tail becomes light, in such case the convergence rate of PE gets slows and NE outperforms PE. On the other hand, NE is demonstrated not to work at all for heavy tailed cases.

Applying the proposed approach, this study considers the following mixture estimator of SMD
\begin{align}
\widehat{G} (x; h) := q G_{\widehat{\bm{\gamma}}}(x) + (1 -q) (\widehat{F}_h)^m(x),
\end{align}
where $q:=q(h)$ is a monotone function satisfying (2). The proposed estimator satisfies
$$\widehat{G} (x; h) \to 
\begin{dcases}
(\bar{F})^m(x) ~~~ & {\rm as} ~~~ h \to 0 \\
G_{\widehat{\bm{\gamma}}}(x) ~~~ & {\rm as} ~~~ h \to \infty.
\end{dcases}
$$

In order to obtain asymptotic convergence rate, we introduce {\rm (i)} (the Hall class) $\alpha>0$, $\beta\ge 2^{-1}$, $A>0$ and $B \neq0$ exist s.t. 
\begin{align*}
x^{\alpha+\beta} \{1 - F(x) - A x^{-\alpha} (1 + B x^{-\beta} ) \} \to 0~~~ {\rm as} ~~~ x \to \infty,
\end{align*}
or {\rm (ii)} (the Weibull class) $\kappa>0$ and $C>0$ exist s.t. 
\begin{align*}
\exp(Cx^{\kappa})\{1 - F(x) - \exp(-Cx^{\kappa})\} \to 0~~~ {\rm as} ~~~ x \to \infty,
\end{align*}
or {\rm (iii)} (the bounded class) $x^* \in \mathbb{R}$, $\mu<-2$, $\sigma\le -2^{-1}$, $D>0$ and $E\neq0$ exist s.t. 
\begin{align*}
(x^* -x)^{\mu+\sigma}\{1 -F(x) - (x^* -x)^{-\mu} (D + E (x^* -x)^{-\sigma})\} \to 0 ~~~ {\rm as} ~~~ x \uparrow x^*.
\end{align*}

Applying the H\"{o}der's inequality to the MSE $\mathbb{E}[\{F^m(x) - \widehat{G} (x; h)\}^2]$ we have the error bound
\begin{align*}
\Psi(x,m,h) := q^2 \mathbb{E}[\{F^m(x) - G_{\widehat{\bm{\gamma}}}(x)\}^2] + (1-q)^2 \mathbb{E}[\{F^m(x) - (\widehat{F}_h)^m(x)\}^2].
\end{align*}
On minimizing the error bound $\Psi(x,m,h)$, we have the following asymptotic result, which is also a consequence of \cite{moriyama2025application}.
\begin{theorem}
Suppose that $(x^* - x)^{-1}h \to 0$ if $h \to 0$. Then, 
\begin{align*}
\Psi(x,m,h) \lesssim 
&\begin{dcases}
q^2 \mathbb{E}[\{F^m(x) -G_{\widehat{\bm{\gamma}}}(x)\}^2] + \mathbb{E}[\{F^m(x) - (\widehat{F}_h)^m(x)\}^2] & {\rm for} ~~~ h \to 0 \\
\mathbb{E}[\{F^m(x) - G_{\widehat{\bm{\gamma}}}(x)\}^2] + (1-q)^2 \mathbb{E}[\{F^m(x) - (\widehat{F}_h)^m(x)\}^2 & {\rm for} ~~~ h \to \infty,
\end{dcases}\\
\intertext{where $\mathbb{E}[\{F^m(x) -\widehat{F}^m(x)\}^2]$ for $h \to 0$ is of order}
&\begin{dcases}
h^{4} m^{-4\gamma} + \frac{m}{n}(1- hm^{-\gamma} ) &{\rm for ~ (i), ~ (iii)}\\
h^{4} (\ln m)^{4\kappa^{-1}(\kappa-1)} + \frac{m}{n} (1- h (\ln m)^{\kappa^{-1}(\kappa-1)} ) &{\rm for ~ (ii)}.
\end{dcases}\\
\intertext{If $\ln F^m(x)$ {$=O(1)$}, $\mathbb{E}[\{F^m(x) - G_{\widehat{\bm{\gamma}}}(x)\}^2]$ is of order}
& \begin{dcases}
m^{2\gamma\rho} + m^{-2} + n^{-1} m(m^{1+2\rho}+1) ~~~ &{\rm for ~ (i), ~ (iii)}\\
C^{-2} (\ln m)^{2} + n^{-1} m(m (\ln m)^{-2} +1) ~~~ &{\rm for ~ (ii)}
\end{dcases}
\end{align*}
{where $\rho$ is the second order parameter}.
\end{theorem}

Suppose $h \to 0$. Then, for {the extreme value index} $\gamma<0$ or $\gamma=0$ with $\kappa>1$, the explicit formula of the optimal bandwidth 
$$h^* := \underset{h} {\operatorname{argmin}} ~ \mathbb{E}[\{F^m(x) -\widehat{F}^m(x)\}^2]$$ 
converging to zero always exists and the optimal order of the MSE is $(m/n)^{2/3}$ (see \citealt{moriyama2025application}). For $\gamma>0$ or $\gamma=0$ with $\kappa\le1$ the explicit formula does not always exist but in such case $m^{-4\gamma} =o(m/n)$ or $h^{4} (\ln m)^{4\kappa^{-1}(\kappa-1)}=o(m/n)$ with any $h \to 0$, that both means $\mathbb{E}[\{F^m(x) -\widehat{F}^m(x)\}^2] = O(m/n)$. 

To summarize under $h \to 0$ the optimal order of the bandwidth $h^*$ cannot always be specified; however, there exist $h$ satisfying $\mathbb{E}[\{F^m(x) -\widehat{F}^m(x)\}^2] = O(m/n)$. The convergence rate is always faster than that of $\mathbb{E}[\{F^m(x) - G_{\widehat{\bm{\gamma}}}(x)\}^2]$ (see also Table {12}), that is finally we have a conclusion
\begin{align*}
\min_{h \to 0} \Psi(x,m,h) \lesssim  \min_{h \to \infty} \Psi(x,m,h).
\end{align*}
Let us consider the rate of the optimal $h$ and suppose $q=h(1+h)^{-1}$ as an simple example. Then, the asymptotic minimizer of $h$ is obtained by
$$\min_h \{h^2 \mathbb{E}[\{F^m(x) -G_{\widehat{\bm{\gamma}}}(x)\}^2] + \mathbb{E}[\{F^m(x) - (\widehat{F}_h)^m(x)\}^2]\}.$$
However, to obtain the explicit form of the solution it is {still} too complicated, and so let us consider another approach. First, $h^*$ is considered to be optimal if $h^*$ exists and
\begin{align}
\frac{n}{m} (h^*)^2 \mathbb{E}[\{F^m(x) -G_{\widehat{\bm{\gamma}}}(x)\}^2]
\end{align}
converges to some constant. Otherwise, the next idea {of finding the optimal bandwidth in some sense} is balancing the two convergence rate
\begin{align}
O(h^2 \mathbb{E}[\{F^m(x) - G_{\widehat{\bm{\gamma}}}(x)\}^2]) = O\left(h^{4} m^{-4\gamma} + \frac{m}{n}\right).
\end{align}

Following the above two-step choice the polynomial degree of the optimal $h$ as a function of $n$ and the convergence rate of the MSE of the mixture SMD estimator in each case {is} summarized in Table 1. The optimal bandwidth value with the asterisk e.g. $-1/3^*$ means that the optimal $h^*$ exists and $h^*=O(n^{-1/3})$. $0-$ means $h^*$ does not exist but 
\begin{align*}
\mathbb{E}[\{F^m(x) - G_{\widehat{\bm{\gamma}}}(x)\}^2] = O\left(m^{-4\gamma} + \frac{m}{n}\right)
\end{align*}
{holds,} that is, the two convergence rate{s of PE and NE are} balanced no matter how slow the convergence rate of $h$ is. $-\infty$ means the rate of the optimal $h$ needs to be enough fast since the rate of $\mathbb{E}[\{F^m(x) -G_{\widehat{\bm{\gamma}}}(x)\}^2]$ cannot be specified and it possibly diverges.

Table 1 shows the convergence rate of the MSE is $m/n$ in all the cases. The rate is not faster than that of the optimal NE, but the better performance in finite sample cases is demonstrated in the next section.

\begin{table}
\caption{The polynomial convergence rates of {an optimal} $h$ and MSE}{\fontsize{6pt}{6pt}\selectfont
$$\begin{tabu}[c]{|ccc||c|c|c||c|c|c|}
 \hline
 \multicolumn{3}{|c||}{{\rm Pareto}} & {m=n^{1/4}} & {m=n^{1/2}} & {m=n^{3/4}} & {m=n^{1/4}} & {m=n^{1/2}} & {m=n^{3/4}} \\ \hline
\ell & \alpha & \beta & \multicolumn{3}{c|| }{{\rm optimal ~} h} & \multicolumn{3}{c|}{{\rm MSE}} \\ \hline

 1/2 & 1/2 & 1 & -1/8 & 0- & 0- & -3/4 & -1/2 & -1/4 \\
 1 & 1 & 1 & -1/8 & 0- & 0- & -3/4 & -1/2 & -1/4 \\
 3 & 3 & 1 & -1/2^* & -1/12 & 0- & -3/4 & -1/2 & -1/4 \\
 10 & 10 & 1 & -27/40^* & -7/20^* & -1/20 & -3/4 & -1/2 & -1/4 \\
 
\hline\hline
 \multicolumn{3}{|c||}{{\rm T}} & {m=n^{1/4}} & {m=n^{1/2}} & {m=n^{3/4}} & {m=n^{1/4}} & {m=n^{1/2}} & {m=n^{3/4}} \\ \hline
\ell & \alpha & \beta & \multicolumn{3}{c|| }{{\rm optimal ~} h} & \multicolumn{3}{c|}{{\rm MSE}} \\ \hline

 1/2 & 1/2 & 1 & -1/8 & 0- & 0- & -3/4 & -1/2 & -1/4 \\
 1 & 1 & 1 & -1/8 & 0- & 0- & -3/4 & -1/2 & -1/4 \\
 3 & 3 & 1 & -5/24 & 0- & 0- & -3/4 & -1/2 & -1/4 \\
 10 & 10 & 1 & -27/40^* & -7/20^* & -1/40^* & -3/4 & -1/2 & -1/4 \\

 \hline\hline
 \multicolumn{3}{|c||}{{\rm Burr}} & {m=n^{1/4}} & {m=n^{1/2}} & {m=n^{3/4}} & {m=n^{1/4}} & {m=n^{1/2}} & {m=n^{3/4}} \\ \hline
c,\ell & \alpha & \beta & \multicolumn{3}{c|| }{{\rm optimal ~} h} & \multicolumn{3}{c|}{{\rm MSE}} \\ \hline

 1/2,1/2 & 1/4 & 1/2 & -1/8 & 0- & 0- & -3/4 & -1/2 & -1/4 \\
 1,1/2 & 1/2 & 1 & -1/8 & 0- & 0- & -3/4 & -1/2 & -1/4 \\
 3,1/2 & 3/2 & 3 & -1/4^* & 0- & 0- & -3/4 & -1/2 & -1/4 \\
 1/2,1 & 1/2 & 1/2 & -1/8 & 0- & 0- & -3/4 & -1/2 & -1/4 \\
 1,1 & 1 & 1 & -1/8 & 0- & 0- & -3/4 & -1/2 & -1/4 \\
 3,1 & 3 & 3 & -1/2^* & 0- & 0- & -3/4 & -1/2 & -1/4 \\
 1/2,3 & 3/2 & 1/2 & -7/24 & -1/12 & 0- & -3/4 & -1/2 & -1/4 \\
 1,3 & 3 & 1 & -1/2^* & -1/12 & 0- & -3/4 & -1/2 & -1/4 \\
 3,3 & 9 & 3 & -2/3^* & -1/3^* & 0- & -3/4 & -1/2 & -1/4 \\
 
 \hline\hline
 \multicolumn{3}{|c||}{{\rm Fr{e}chet}} & {m=n^{1/4}} & {m=n^{1/2}} & {m=n^{3/4}} & {m=n^{1/4}} & {m=n^{1/2}} & {m=n^{3/4}} \\ \hline
\gamma & \alpha & \beta & \multicolumn{3}{c|| }{{\rm optimal ~} h} & \multicolumn{3}{c|}{{\rm MSE}} \\ \hline

 5 & 1/5 & 1/5 & -\infty & -\infty & -\infty & -3/4 & -1/2 & -1/4 \\
 2 & 1/2 & 1/2 & -1/8 & 0- & 0- & -3/4 & -1/2 & -1/4 \\
 1 & 1 & 1 & -1/8 & 0- & 0- & -3/4 & -1/2 & -1/4 \\
 1/2 & 2 & 1 & -3/8^* & 0- & 0- & -3/4 & -1/2 & -1/4 \\
 1/4 & 4 & 1 & -9/16^* & -1/4 & 0- & -3/4 & -1/2 & -1/4 \\
 
 \hline\hline
 \multicolumn{3}{|c||}{{\rm Weibull}} & {m=n^{1/4}} & {m=n^{1/2}} & {m=n^{3/4}} & {m=n^{1/4}} & {m=n^{1/2}} & {m=n^{3/4}} \\ \hline
\kappa& \gamma & \rho & \multicolumn{3}{c|| }{{\rm optimal ~} h} & \multicolumn{3}{c|}{{\rm MSE}} \\ \hline

1/2 & 0 & 0 & -\infty & -\infty & -\infty & -3/4 & -1/2 & -1/4 \\
1 & 0 & 0 & -\infty & -\infty & -\infty & -3/4 & -1/2 & -1/4 \\
3 & 0 & 0 & -\infty & -\infty & -\infty & -3/4 & -1/2 & -1/4 \\
10 & 0 & 0 & -\infty & -\infty & -\infty & -3/4 & -1/2 & -1/4 \\

 \hline\hline
 \multicolumn{3}{|c||}{{\rm inv. Burr}} & {m=n^{1/4}} & {m=n^{1/2}} & {m=n^{3/4}} & {m=n^{1/4}} & {m=n^{1/2}} & {m=n^{3/4}} \\ \hline
c,\ell & \mu & \sigma & \multicolumn{3}{c|| }{{\rm optimal ~} h} & \multicolumn{3}{c|}{{\rm MSE}} \\ \hline

3, 2 & -6 & -2 & -7/8^* & -3/4^* & -5/8^* & -3/4 & -1/2 & -1/4 \\
1, 2 & -2 & -2 & -\infty & -\infty & -\infty & -3/4 & -1/2 & -1/4 \\
 1/2, 2 & -1 & -2 & -\infty & -\infty & -\infty & -3/4 & -1/2 & -1/4 \\
 3, 1 & -3 & -1 & -1^* & -1^* & -1^* & -3/4 & -1/2 & -1/4 \\
 1, 1 & -1 & -1 & -\infty & -\infty & -\infty & -3/4 & -1/2 & -1/4 \\

 \hline
\end{tabu}$$
}
\end{table}

\section{Simulation study on SMD estimation}

This section investigates numerical properties of the proposed mixture through simulation experiments in SMD estimation. The naive mixture approach proposed by \cite{schuster1985parametric} and \cite{olkin1987semiparametric} constructs the following estimator
$$\widetilde{G} (x; p, h) := p G_{\widehat{\bm{\gamma}}}(x) + (1- p) (\widehat{F}_{\widehat{h}}^m)(x),$$
where the mixing ratio $p$ is replaced with
$$\widehat{p} :=\mathop{\rm argmax}\limits_{0\le p\le1} \left[ \sum_{j=1}^N \log \widetilde{g} (Y_j; p, \widehat{h}) \right].$$
$\widetilde{g}$ is the partial derivative of $\widetilde{G}$ with respect to $x$. $\widehat{h}$ needs to be a data-driven bandwidth determined before the maximization. It has been proven that the optimal bandwidth of the pointwise mean-squared error (MSE) for $\widehat{F}_h$ is also optimal for {NE in some sense} (\citealt{moriyama2025application}). This study employs the estimators proposed in \cite{altman1995bandwidth} or \cite{bowman1998bandwidth}. The \texttt{kerdiest} package in R provides the \texttt{ALbw} function and the \texttt{CVbw} function for calculating {each} estimated value. 

The proposed mixture estimator needs to determine the tuning parameter, which in this study section is aimed at minimizing the integrated squared error
$$\int \{\widehat{G}(x; h) -F^m(x)\}^2 {\rm d}x$$
{or minimizing Anderson-Darling-type metric given by}
$$\int \frac{\{\widehat{G}(x;h) - F^m(x)\}^2}{F^m(x)\{1-F^m(x)\}} {\rm d}F^m(x),$$
{which is a weighted mean integrated squared error (MISE). The weight function $[F^m(x)\{1-F^m(x)\}]^{-1}$ magnifies the error in the distributional tails.}

This study considers the following cross-validation approach{es}
$$\widehat{h} :=\mathop{\rm argmin}\limits_{h>0} \left[\frac{1}{N} \sum_{j=1}^N \int \{ \widehat{G}^{(-j)}(x; h) - I(x>Y_j) \}^2 {\rm d}x \right]$$
{or}
$$\widehat{h} :=\mathop{\rm argmin}\limits_{h>0} \left[\frac{1}{N} \sum_{j=1}^N \frac{\{ (N+1)\widehat{G}^{(-j)}(Y_{(j)}; h) - j \}^2}{j(N+1-j)} \right],$$
{based on the order statistic of the block maxima $\{Y_{(1)} < \cdots < Y_{(N)}\}$, where $\widehat{G}^{(-j)}$ is the estimated mixture distribution without $\{X_{m(j-1)+1}, \cdots, X_{mj}\}$ and $I$ is the indicator function.}

To simulate the following MISE 
$$L_m^{-1} \int_{Q_ m(0.1)}^{Q_ m(0.9)} \left(\bar{G}(x) -F^m(x)\right)^2 {\rm d}x,$$
we surveyed the numerical accuracy of the two mixture estimators in finite-sample cases where $L_m :=Q_ m(0.9) - Q_ m(0.1)$ and $Q_ m(r)$ denotes the $r$th quantile of the SMD. $\bar{G}$ is the naive mixture $\widetilde{G}$ or the proposed $\widehat{G}$.

We simulated MISE values 100 times between the $10$th and $90$th quantiles of the SMD. Tables {2--5} denote the mean MISE values and standard deviation (sd), respectively. The sample sizes were $(n=)2^{8}$ or $2^{12}$, and $m =n^{1/4}$, $m=n^{1/2}$, and $m=n^{3/4}$. The underlying distributions, $F$, were Pareto, T, Burr, Frechet, Weibull, inverse Burr, and non-MDA {given by} 
$$1-F(x)=e^{-x -\sin x} ~~~x>0,$$ 
which does not belong to the maximum domain of attraction ({non-MDA}, {\citealt{dehaan2006extreme}}). Kernels $w$ {were} of the Epanechnikov type for inverse Burr distributions and Gaussian for the others. The mixing ratio $p$ of the naive mixture $\widetilde{G}$ are provided in Table{s 6--7} and $q=h(1+h)^{-1}$ of the proposed $\widehat{G}$ in Table{s 8--9}. In the columns named P/N in Table{s 2--3}, P means that the mean of the estimated mixing ratio values $p$ was greater than 0.7 (i.e. the mixture estimator being close to parametric). Conversely, N means the value was less than 0.3. The blank means the mean value is between 0.3 and 0.7. Table{s 4--5} shows that of the mixing ratio $q=h(1+h)^{-1}$.

The parametric distribution estimator {was} demonstrated to be better than the nonparametric distribution estimator in heavy tailed cases (Tables {13--16} in the supplementary file, which was provided by \citealt{moriyama2025application}). The columns P/N show that {the naive mixture estimator} generally becomes parametric as the tail gets heavy. The numerical property {coincided} with what we expected in advance. {For the cases with $\gamma\le0$ the naive mixture estimator was close to the nonparametric one in general. The difference between ALbw and CVbw was not so large, but CVbw-based mixture estimator was more numerically stable for heavy-tailed cases.}

{The proposed estimator minimizing the integrated squared error (ISE) behaved differently depending on the extreme value index $\gamma$. ISE-based estimator tended to be parametric for $\gamma\le0$ and nonparametric for $\gamma>0$.} However, {in such cases} the proposed estimator does not necessarily behave like the wholly nonparametric one since the mixture estimator with $q = 0$ does not coincide with the naive kernel estimator. Though ISE-based estimator was almost nonparametric for heavy-tailed cases, the numerical performance is much better than the naive kernel estimator (see also Tables 13--15), and ISE-based estimator outperformed even the naive mixture estimators in some times.

The Anderson-Darling-type estimator (AD) became rather close to PE in general, which happened even when the numerical superiority of the nonparametric estimator in such cases is demonstrated. The results are not preferable, especially for $\gamma \fallingdotseq 0$ and the non-MDA case as seen by comparing Tables 2--5. As far as $\gamma=0$, the naive mixture estimator especially ALbw-based estimator seems to perform best among them. For the non-MDA case the mixture estimators except AD-based estimator became nonparametric.

It was confirmed AD-based estimator sometimes outperforms ISE-based estimator in the sense of MISE. The difference between the results of ISE-based and AD-based estimators is considered to come from the weight function. The result suggests the PE captures the tail behavior of SMD more than NE in such cases. 

By examining Tables 2--5, we can summarize the results as follows. For heavy-tailed cases or $\gamma<0$ but close to zero the proposed estimators are recommended. When $\gamma = 0$ the naive mixture estimators especially ALbw-based estimator is considered to be acceptable. If $F$ seems not to belong to any MDA, AD-based estimator should be avoided.

Numerical experiments showed that the computation time of {ALbw-based, CVbw-based, AD-based estimators took  around 0.1, 100, 0.7 times than that of ISE-based one respectively for $n=2^{8}$. These were around 10, 2000, 2 times for $n=2^{12}$}. \texttt{ALbw} function is computationally very efficient method indeed, but the calculation cost of the bandwidth selection rapidly increases as $n$ gets large in general. The proposed mixture approach succeeded in substantial reduction of the computational cost. From the above this section concludes that the proposed estimators for SMD are demonstrated to be enough effective.

\begin{landscape}
\begin{table}
\caption{{Scaled MISE values ($\times 100$) and sd values ($\times 100$) of the naive mixture approach with ALbw function}}
{\fontsize{6pt}{6pt}\selectfont
$$\begin{tabu}[c]{|c||ccc|ccc|ccc||ccc|ccc|ccc|}
 \hline
 & \multicolumn{9}{c||}{n=2^8} & \multicolumn{9}{c|}{n=2^{12}} \\ \hline
 & {\rm MISE} & {\rm sd} & {\rm P/N} & {\rm MISE} & {\rm sd} & {\rm P/N} & {\rm MISE} & {\rm sd} & {\rm P/N} & {\rm MISE} & {\rm sd} & {\rm P/N} & {\rm MISE} & {\rm sd} & {\rm P/N} & {\rm MISE} & {\rm sd} & {\rm P/N} \\ \hline\hline 
 {\rm Pareto} & \multicolumn{3}{c|}{m=2^{2}} & \multicolumn{3}{c|}{m=2^{4}} & \multicolumn{3}{c||}{m=2^{6}} & \multicolumn{3}{c|}{m=2^{3}} & \multicolumn{3}{c|}{m=2^{6}} & \multicolumn{3}{c|}{m=2^{9}} \\ \hline
1/2 & 0.315 & 1.498 & {\rm P} & 4.550 & 6.869 & {\rm P} & 11.94 & 12.34 & {\rm P} & 0.276 & 1.420 & {\rm P} & 0.562 & 2.386 & {\rm P} & 8.426 & 9.304 & {\rm P} \\
1 & 0.282 & 0.447 & & 1.085 & 2.379 & & 8.187 & 6.003 & & 0.346 & 0.883 & & 0.264 & 0.888 & & 5.646 & 5.586 & \\
3 & 0.195 & 0.200 & {\rm N} & 0.813 & 0.808 & & 4.502 & 3.940 & & 0.027 & 0.026 & {\rm N} & 0.199 & 0.213 & {\rm N} & 1.923 & 1.838 & \\
10 & 0.201 & 0.217 & {\rm N} & 0.839 & 0.835 & {\rm N} & 4.501 & 4.010 & & 0.026 & 0.025 & {\rm N} & 0.211 & 0.219 & {\rm N} & 1.768 & 1.689 & \\

 \hline
 {\rm T} & \multicolumn{3}{c|}{m=2^{2}} & \multicolumn{3}{c|}{m=2^{4}} & \multicolumn{3}{c||}{m=2^{6}} & \multicolumn{3}{c|}{m=2^{3}} & \multicolumn{3}{c|}{m=2^{6}} & \multicolumn{3}{c|}{m=2^{9}}\\ \hline
1/2 & 0.321 & 1.691 & {\rm P} & 4.170 & 6.881 & {\rm P} & 12.48 & 12.46 & {\rm P} & 0.109 & 0.880 & {\rm P} & 0.395 & 1.797 & & 9.310 & 9.739 & {\rm P}\\
1 & 0.243 & 0.282 & {\rm P} & 0.732 & 1.026 & {\rm P} & 7.807 & 6.381 & & 0.029 & 0.032 & {\rm P} & 0.166 & 0.200 & {\rm P} & 5.556 & 5.687 & {\rm P}\\
3 & 0.200 & 0.209 & & 0.823 & 0.860 & & 4.434 & 3.942 & & 0.032 & 0.034 & {\rm P} & 0.188 & 0.218 & {\rm P} & 1.865 & 1.719 & \\
10 & 0.189 & 0.215 & {\rm N} & 0.745 & 0.827 & {\rm N} & 4.032 & 4.166 & & 0.023 & 0.025 & {\rm N} & 0.182 & 0.206 & {\rm N} & 1.740 & 1.711 & \\

 \hline
 {\rm Burr} & \multicolumn{3}{c|}{m=2^{2}} & \multicolumn{3}{c|}{m=2^{4}} & \multicolumn{3}{c||}{m=2^{6}} & \multicolumn{3}{c|}{m=2^{3}} & \multicolumn{3}{c|}{m=2^{6}} & \multicolumn{3}{c|}{m=2^{9}}\\ \hline
1/2,1/2 & 0.710 & 2.988 & {\rm P} & 4.192 & 9.591 & {\rm P} & 15.66 & 20.20 & {\rm P} & 3.488 & 8.120 & {\rm P} & 1.037 & 1.742 & {\rm P} & 0.382 & 1.278 & \\ 
 1,1/2 & 0.266 & 1.262 & {\rm P} & 4.550 & 7.384 & {\rm P} & 12.39 & 13.33 & {\rm P} & 0.262 & 0.627 & {\rm P} & 0.405 & 1.055 & {\rm P} & 4.102 & 6.866 & \\ 
 3,1/2 & 0.208 & 0.244 & & 0.836 & 1.097 & & 6.333 & 4.532 & & 0.030 & 0.172 & & 0.199 & 0.352 & & 2.745 & 3.517 & \\ 
 1/2,1 & 0.339 & 1.581 & {\rm P} & 4.608 & 6.733 & {\rm P} & 12.44 & 12.92 & {\rm P} & 0.323 & 0.775 & {\rm P} & 0.484 & 1.422 & {\rm P} & 4.567 & 6.373 & \\ 
 1,1 & 0.271 & 0.386 & & 1.277 & 3.004 & & 8.226 & 6.478 & & 0.307 & 0.824 & & 0.444 & 1.692 & & 4.310 & 4.974 & \\ 
 3,1 & 0.177 & 0.192 & {\rm N} & 0.776 & 0.773 & & 4.645 & 3.846 & & 0.026 & 0.026 & {\rm N} & 0.201 & 0.213 & & 1.931 & 1.975 & \\ 
 1/2,3 & 0.200 & 0.305 & {\rm P} & 0.887 & 1.314 & & 6.432 & 4.798 & & 0.027 & 0.037 & & 0.196 & 0.245 & & 2.972 & 3.732 & \\ 
 1,3 & 0.204 & 0.214 & {\rm N} & 0.851 & 0.858 & & 4.880 & 4.080 & & 0.026 & 0.026 & {\rm N} & 0.200 & 0.189 & {\rm N} & 1.864 & 1.957 & \\ 
 3,3 & 0.188 & 0.209 & {\rm N} & 0.729 & 0.775 & {\rm N} & 3.915 & 3.904 & & 0.024 & 0.025 & {\rm N} & 0.198 & 0.200 & {\rm N} & 1.740 & 1.858 & \\ 

 \hline
 {\rm Fr{e}chet} & \multicolumn{3}{c|}{m=2^{2}} & \multicolumn{3}{c|}{m=2^{4}} & \multicolumn{3}{c||}{m=2^{6}} & \multicolumn{3}{c|}{m=2^{3}} & \multicolumn{3}{c|}{m=2^{6}} & \multicolumn{3}{c|}{m=2^{9}} \\ \hline
5 & 1.398 & 6.775 & {\rm P} & 3.676 & 9.607 & {\rm P} & 16.80 & 22.75 & {\rm P} & 11.54 & 17.14 & {\rm P} & 2.470 & 7.609 & {\rm P} & 14.69 & 22.68 & {\rm P} \\
2 & 0.342 & 1.914 & {\rm P} & 4.406 & 7.427 & {\rm P} & 12.79 & 13.21 & {\rm P} & 0.242 & 1.195 & {\rm P} & 0.445 & 1.979 & {\rm P} & 9.449 & 10.36 & {\rm P} \\
1 & 0.256 & 0.372 & {\rm P} & 1.262 & 2.729 & & 8.386 & 6.756 & & 0.205 & 0.573 & & 0.262 & 0.923 & & 5.619 & 5.603 & \\
1/2 & 0.212 & 0.240 & {\rm N} & 0.844 & 0.815 & & 5.621 & 4.226 & & 0.026 & 0.026 & {\rm N} & 0.189 & 0.199 & & 2.206 & 2.489 & \\
1/4 & 0.191 & 0.210 & {\rm N} & 0.775 & 0.755 & & 4.710 & 4.286 & & 0.026 & 0.028 & {\rm N} & 0.207 & 0.212 & {\rm N} & 1.823 & 1.708 & \\
 
 \hline
 {\rm Weibull} & \multicolumn{3}{c|}{m=2^{2}} & \multicolumn{3}{c|}{m=2^{4}} & \multicolumn{3}{c||}{m=2^{6}} & \multicolumn{3}{c|}{m=2^{3}} & \multicolumn{3}{c|}{m=2^{6}} & \multicolumn{3}{c|}{m=2^{9}}\\ \hline
1/2 & 0.195 & 0.218 & & 0.849 & 0.842 & {\rm N} & 4.910 & 4.180 & & 0.025 & 0.023 & {\rm N} & 0.203 & 0.201 & {\rm N} & 1.747 & 1.598 & \\
1 & 0.208 & 0.219 & {\rm N} & 0.853 & 0.900 & {\rm N} & 4.636 & 4.570 & & 0.026 & 0.026 & {\rm N} & 0.191 & 0.198 & {\rm N} & 1.754 & 1.592 & \\
3 & 0.179 & 0.211 & {\rm N} & 0.680 & 0.860 & {\rm N} & 3.455 & 4.034 & & 0.023 & 0.024 & {\rm N} & 0.168 & 0.192 & {\rm N} & 1.597 & 1.727 & \\
10 & 0.175 & 0.203 & {\rm N} & 0.732 & 0.863 & {\rm N} & 3.567 & 4.218 & & 0.026 & 0.029 & {\rm N} & 0.190 & 0.224 & {\rm N} & 1.515 & 1.669 & \\
 
\hline
 {\rm inv. Burr} & \multicolumn{3}{c|}{m=2^{2}} & \multicolumn{3}{c|}{m=2^{4}} & \multicolumn{3}{c||}{m=2^{6}} & \multicolumn{3}{c|}{m=2^{3}} & \multicolumn{3}{c|}{m=2^{6}} & \multicolumn{3}{c|}{m=2^{9}}\\ \hline
 3,2 & 12.00 & 2.409 & {\rm N} & 27.14 & 5.582 & {\rm N} & 21.88 & 12.54 & & 20.34 & 0.940 & {\rm N} & 31.88 & 0.214 & {\rm N} & 25.30 & 13.66 & {\rm N} \\
1,2 & 11.77 & 2.159 & {\rm N} & 25.50 & 3.147 & {\rm N} & 20.53 & 9.972 & & 19.38 & 1.190 & {\rm N} & 27.28 & 0.497 & {\rm N} & 25.28 & 6.694 & {\rm N} \\
1/2,2 & 10.31 & 1.524 & {\rm N} & 18.66 & 1.475 & {\rm N} & 14.99 & 6.288 & {\rm N} & 15.28 & 2.976 & {\rm N} & 18.66 & 0.558 & {\rm N} & 17.50 & 2.848 & {\rm N} \\
3,1 & 5.789 & 3.558 & & 12.74 & 6.202 & & 10.35 & 8.956 & & 11.54 & 5.843 & & 17.94 & 3.946 & {\rm N} & 13.04 & 10.71 & \\
 1,1 & 4.204 & 2.798 & & 9.215 & 4.600 & & 10.23 & 7.253 & & 7.440 & 3.970 & & 9.700 & 3.485 & {\rm N} & 8.891 & 6.122 & \\
  1/2,1 & 2.647 & 2.376 & & 4.417 & 2.630 & & 6.682 & 4.234 & {\rm N} & 3.725 & 2.697 & & 4.870 & 2.022 & & 7.045 & 2.548 & {\rm N} \\
 3,1/3 & 6.655 & 5.355 & & 3.070 & 2.678 & & 8.522 & 4.816 & & 2.178 & 0.887 & & 6.556 & 2.210 & & 10.38 & 4.666 & {\rm N} \\
 1,1/3 & 11.60 & 7.117 & & 2.300 & 2.807 & & 4.351 & 3.969 & {\rm N} & 1.297 & 1.203 & & 5.188 & 4.181 & & 4.628 & 1.616 & {\rm N} \\
 1/2,1/3 & 13.23 & 10.21 & & 3.132 & 5.099 & & 2.481 & 1.683 & {\rm N} & 4.467 & 5.464 & & 8.356 & 4.243 & & 2.743 & 2.466 & {\rm N} \\
 
\hline
 
 \text{non-MDA} & \multicolumn{3}{c|}{m=2^{2}} & \multicolumn{3}{c|}{m=2^{4}} & \multicolumn{3}{c||}{m=2^{6}} & \multicolumn{3}{c|}{m=2^{3}} & \multicolumn{3}{c|}{m=2^{6}} & \multicolumn{3}{c|}{m=2^{9}}\\ \hline
 & 0.194 & 0.201 & {\rm N} & 1.591 & 1.769 & {\rm N} & 3.789 & 4.257 & & 0.041 & 0.096 & {\rm N} & 0.183 & 0.180 & {\rm N} & 1.959 & 1.855 & \\

\hline
 \end{tabu}$$
}
\end{table}
\end{landscape}

\clearpage

\begin{landscape}
\begin{table}
\caption{{Scaled MISE values ($\times 100$) and sd values ($\times 100$) of the naive mixture approach with CVbw function}}
{\fontsize{6pt}{6pt}\selectfont
$$\begin{tabu}[c]{|c||ccc|ccc|ccc||ccc|ccc|ccc|}
 \hline
 & \multicolumn{9}{c||}{n=2^8} & \multicolumn{9}{c|}{n=2^{12}} \\ \hline
 & {\rm MISE} & {\rm sd} & {\rm P/N} & {\rm MISE} & {\rm sd} & {\rm P/N} & {\rm MISE} & {\rm sd} & {\rm P/N} & {\rm MISE} & {\rm sd} & {\rm P/N} & {\rm MISE} & {\rm sd} & {\rm P/N} & {\rm MISE} & {\rm sd} & {\rm P/N} \\ \hline\hline 
 {\rm Pareto} & \multicolumn{3}{c|}{m=2^{2}} & \multicolumn{3}{c|}{m=2^{4}} & \multicolumn{3}{c||}{m=2^{6}} & \multicolumn{3}{c|}{m=2^{3}} & \multicolumn{3}{c|}{m=2^{6}} & \multicolumn{3}{c|}{m=2^{9}} \\ \hline
1/2 & 0.237 & 1.078 & {\rm P} & 4.598 & 7.112 & {\rm P} & 12.04 & 12.79 & {\rm P} & 0.354 & 1.111 & {\rm P} & 0.688 & 2.397 & {\rm P} & 3.926 & 6.411 & \\
  1 & 0.245 & 0.321 &  & 1.106 & 2.377 &  & 8.259 & 6.582 & {\rm P} & 2.056 & 1.941 & & 0.781 & 1.770 & & 3.107 & 4.699 & \\
  3 & 0.212 & 0.224 & {\rm N} & 0.836 & 0.833 &  & 4.624 & 3.787 &  & 0.040 & 0.083 & {\rm N} & 0.233 & 0.251 & {\rm N} & 1.623 & 1.433 & \\
  10 & 0.212 & 0.222 & {\rm N} & 0.850 & 0.907 & {\rm N} & 4.361 & 3.790 &  & 0.026 & 0.026 & {\rm N} & 0.221 & 0.192 & {\rm N} & 1.815 & 1.867 & {\rm N} \\ 

 \hline
 {\rm T} & \multicolumn{3}{c|}{m=2^{2}} & \multicolumn{3}{c|}{m=2^{4}} & \multicolumn{3}{c||}{m=2^{6}} & \multicolumn{3}{c|}{m=2^{3}} & \multicolumn{3}{c|}{m=2^{6}} & \multicolumn{3}{c|}{m=2^{9}}\\ \hline
1/2 & 0.213 & 0.225 & {\rm P} & 3.709 & 6.819 & {\rm P} & 12.89 & 13.11 & {\rm P} & 0.216 & 0.402 & {\rm P} & 0.459 & 1.562 & {\rm P} & 4.279 & 6.860 &  \\ 
  1 & 0.243 & 0.276 & {\rm P} & 0.863 & 1.827 & {\rm P} & 7.775 & 6.077 & & 0.044 & 0.158 & {\rm P} & 0.178 & 0.199 & {\rm P} & 4.727 & 5.153 & \\ 
  3 & 0.196 & 0.227 &  & 0.805 & 0.800 &  & 4.668 & 4.173 & & 0.037 & 0.033 & {\rm P} & 0.228 & 0.212 & {\rm P} & 1.528 & 1.476 & \\ 
  10 & 0.188 & 0.215 & {\rm N} & 0.760 & 0.843 & {\rm N} & 3.955 & 4.280 & & 0.026 & 0.026 & {\rm N} & 0.184 & 0.184 & {\rm N} & 1.534 & 1.588 & \\ 

 \hline
 {\rm Burr} & \multicolumn{3}{c|}{m=2^{2}} & \multicolumn{3}{c|}{m=2^{4}} & \multicolumn{3}{c||}{m=2^{6}} & \multicolumn{3}{c|}{m=2^{3}} & \multicolumn{3}{c|}{m=2^{6}} & \multicolumn{3}{c|}{m=2^{9}}\\ \hline
1/2,1/2 & 1.022 & 2.374 & {\rm P} & 1.322 & 4.151 & & 1.501 & 6.026 & {\rm N} & 2.633 & 3.959 & {\rm P} & 1.000 & 0.843 & {\rm P} & 0.242 & 0.671 & {\rm N} \\ 
  1,1/2 & 0.419 & 1.506 & {\rm P} & 2.961 & 5.477 & {\rm P} & 3.183 & 8.251 & {\rm N} & 0.304 & 0.614 & {\rm P} & 0.412 & 1.457 & {\rm P} & 4.548 & 9.089 & \\ 
  3,1/2 & 0.209 & 0.248 & & 0.817 & 1.007 & & 6.053 & 5.909 & & 0.092 & 0.137 & {\rm P} & 0.564 & 1.348 & {\rm N} & 3.037 &  3.204 & \\ 
  1/2,1 & 0.466 & 1.475 & {\rm P} & 3.065 & 5.152 & {\rm P} & 6.023 & 7.230 & & 0.205 & 0.691 & {\rm P} & 0.353 & 1.255 & {\rm P} & 4.458 & 6.810 & \\
  1,1 & 0.247 & 0.353 & & 1.503 & 3.162 & & 8.104 & 5.957 & {\rm P} & 1.763 & 1.725 & & 0.729 & 2.015 & & 3.399 & 4.690 & \\ 
  3,1 & 0.208 & 0.222 & {\rm N} & 0.872 & 0.887 & & 5.729 & 5.035 & & 0.087 & 0.137 & {\rm P} & 0.331 & 0.595 & & 3.022 & 3.805 & \\ 
  1/2,3 & 0.179 & 0.220 & {\rm P} & 0.854 & 1.175 & & 6.905 & 5.448 & {\rm P} & 0.261 & 0.550 & & 0.261 & 0.531 & & 2.954 & 3.781 & \\
  1,3 & 0.210 & 0.222 & {\rm N} & 0.811 & 0.772 & & 5.459 & 4.770 & {\rm P} & 0.026 & 0.024 & {\rm N} & 0.192 & 0.184 & & 1.730 & 1.565 & \\ 
  3,3 & 0.167 & 0.179 & {\rm N} & 0.760 & 0.811 & {\rm N} & 4.991 & 4.713 & & 0.028 & 0.026 & {\rm N} & 0.185 & 0.198 & {\rm N} & 1.398 & 1.634 & \\ 

 \hline
 {\rm Fr{e}chet} & \multicolumn{3}{c|}{m=2^{2}} & \multicolumn{3}{c|}{m=2^{4}} & \multicolumn{3}{c||}{m=2^{6}} & \multicolumn{3}{c|}{m=2^{3}} & \multicolumn{3}{c|}{m=2^{6}} & \multicolumn{3}{c|}{m=2^{9}} \\ \hline
   5 & 1.251 & 3.309 & {\rm P} & 2.231 & 5.218 & {\rm P} & 6.043 & 6.615 & & 4.801 & 10.78 & {\rm P} & 0.947 & 0.658 & {\rm P} & 0.164 & 0.403 & {\rm N} \\
   2 & 0.443 & 1.687 & {\rm P} & 3.054 & 5.138 & {\rm P} & 7.186 & 6.640 & & 0.225 & 0.451 & {\rm P} & 0.451 & 1.198 & {\rm P} & 3.065 & 5.925 & \\ 
  1 & 0.224 & 0.301 & & 1.345 & 2.791 & & 7.715 & 5.244 & {\rm P} & 1.467 & 1.636 & {\rm P} & 0.753 & 2.071 & & 3.506 & 4.998 & \\ 
  1/2 & 0.217 & 0.280 & {\rm N} & 0.780 & 0.783 & & 5.587 & 4.657 & {\rm P} & 0.186 & 0.588 & {\rm N} & 0.173 & 0.166 & {\rm N} & 2.017 & 2.606 & \\
 1/4 & 0.196 & 0.203 & {\rm N} & 0.804 & 0.865 & & 3.781 & 4.403 & & 0.028 & 0.030 & {\rm N} & 0.177 & 0.159 & & 1.670 & 1.509 & \\
 
  \hline
 {\rm Weibull} & \multicolumn{3}{c|}{m=2^{2}} & \multicolumn{3}{c|}{m=2^{4}} & \multicolumn{3}{c||}{m=2^{6}} & \multicolumn{3}{c|}{m=2^{3}} & \multicolumn{3}{c|}{m=2^{6}} & \multicolumn{3}{c|}{m=2^{9}}\\ \hline
1/2 & 1.189 & 0.648 & & 0.990 & 1.120 & {\rm N} & 5.303 & 3.832 & & 0.495 & 0.232 & & 0.204 & 0.200 & {\rm N} & 15.85 & 3.986 & {\rm P} \\
  1 & 0.567 & 0.591 & {\rm N} & 0.923 & 1.163 & {\rm N} & 4.734 & 4.114 & & 0.025 & 0.028 & {\rm N} & 0.204 & 0.221 & {\rm N} & 1.974 & 1.654 & \\
  3 & 0.323 & 0.373 & {\rm N} & 0.683 & 0.846 & {\rm N} & 3.434 & 3.643 & & 0.025 & 0.025 & {\rm N} & 0.171 & 0.193 & {\rm N} & 9.343 & 5.890 & \\ 
  10 & 0.324 & 0.341 & {\rm N} & 0.756 & 1.013 & {\rm N} & 2.988 & 3.266 & & 0.026 & 0.025 & {\rm N} & 0.183 & 0.161 & {\rm N} & 5.468 & 4.784 & \\
 
\hline
 {\rm inv. Burr} & \multicolumn{3}{c|}{m=2^{2}} & \multicolumn{3}{c|}{m=2^{4}} & \multicolumn{3}{c||}{m=2^{6}} & \multicolumn{3}{c|}{m=2^{3}} & \multicolumn{3}{c|}{m=2^{6}} & \multicolumn{3}{c|}{m=2^{9}}\\ \hline 
  3,2 & 0.363 & 0.647 & {\rm N} & 3.038 & 4.576 & {\rm N} & 7.009 & 7.879 & & 6.457 & 4.173 & {\rm N} & 24.05 & 5.235 & {\rm N} & 13.53 & 13.99 & \\
      1,2 & 0.308 & 0.514 & {\rm N} & 3.035 & 3.973 &  {\rm N} & 8.088 & 7.122 & & 5.049 & 3.475 & {\rm N} & 21.05 & 4.385 & {\rm N} & 10.68 & 12.35 & {\rm N} \\ 
        1/2,2 & 1.465 & 4.625 & {\rm N} & 3.121 & 3.536 & {\rm N} & 7.123 & 6.105 & & 3.654 & 2.939 & {\rm N} & 6.928 & 7.385 & {\rm N} & 4.884 & 7.265 & {\rm N} \\ 
  3,1 & 1.717 & 1.903 &  & 10.40 & 6.750 &  {\rm N} & 7.230 & 8.974 & & 4.275 & 1.994 & & 14.15 & 3.344 & {\rm N} & 10.00 & 10.71 & \\ 
      1,1 & 2.516 & 4.978 & {\rm N} & 5.693 & 5.145 & {\rm N} & 7.898 & 7.316 & & 0.977 & 1.698 & {\rm N} & 3.432 & 4.261 & {\rm N} & 3.398 & 6.446 & {\rm N} \\ 
  1/2,1 & 3.756 & 6.871 & & 3.428 & 4.036 & {\rm N} & 4.592 & 4.814 & {\rm N} & 1.273 & 2.248 & {\rm N} & 2.656 & 2.607 & & 3.994 & 3.969 & {\rm N} \\ 
3,1/3 & 1.464 & 4.589 & {\rm N} & 0.832 & 2.037 &  {\rm N} & 2.865 & 5.256 & {\rm N} & 0.015 & 0.146 & {\rm N} & 2.487 & 3.473 & {\rm N} & 3.763 & 5.926 & {\rm N} \\ 
   1,1/3 & 2.692 & 5.507 & {\rm N} & 2.917 & 4.202 & {\rm N} & 5.935 & 6.385 & & 0.000 & 0.000 & {\rm N} & 4.459 & 2.924 & & 3.860 & 2.070 & {\rm N} \\ 
  1/2,1/3 & 4.866 & 7.591 & & 3.085 & 4.693 & & 2.580 & 2.134 & {\rm N} & 1.416 & 3.600 & {\rm N} & 5.075 & 5.519 & {\rm P} & 2.465 & 1.655 & {\rm N} \\ 
\hline
 
 \text{non-MDA} & \multicolumn{3}{c|}{m=2^{2}} & \multicolumn{3}{c|}{m=2^{4}} & \multicolumn{3}{c||}{m=2^{6}} & \multicolumn{3}{c|}{m=2^{3}} & \multicolumn{3}{c|}{m=2^{6}} & \multicolumn{3}{c|}{m=2^{9}}\\ \hline
 & 0.208 & 0.229 & {\rm N} & 1.528 & 1.761 & {\rm N} & 2.653 & 4.329 & & 0.045 & 0.104 & {\rm N} & 0.182 & 0.178 & {\rm N} & 1.918 & 1.828 & \\

\hline
 \end{tabu}$$
}
\end{table}
\end{landscape}

\clearpage

\begin{landscape}
\begin{table}
\caption{{Scaled MISE values ($\times 100$) and sd values ($\times 100$) of the proposed approach with ISE metric}}
{\fontsize{6pt}{6pt}\selectfont
$$\begin{tabu}[c]{|c||ccc|ccc|ccc||ccc|ccc|ccc|}
 \hline
 & \multicolumn{9}{c||}{n=2^8} & \multicolumn{9}{c|}{n=2^{12}} \\ \hline
 & {\rm MISE} & {\rm sd} & {\rm P/N} & {\rm MISE} & {\rm sd} & {\rm P/N} & {\rm MISE} & {\rm sd} & {\rm P/N} & {\rm MISE} & {\rm sd} & {\rm P/N} & {\rm MISE} & {\rm sd} & {\rm P/N} & {\rm MISE} & {\rm sd} & {\rm P/N} \\ \hline\hline 
 {\rm Pareto} & \multicolumn{3}{c|}{m=2^{2}} & \multicolumn{3}{c|}{m=2^{4}} & \multicolumn{3}{c||}{m=2^{6}} & \multicolumn{3}{c|}{m=2^{3}} & \multicolumn{3}{c|}{m=2^{6}} & \multicolumn{3}{c|}{m=2^{9}} \\ \hline
1/2 & 0.177 & 0.200 & {\rm N} & 1.074 & 2.603 & {\rm N} & 1.734 &  4.779 & {\rm N} & 0.033 & 0.234 & {\rm N} & 0.204 & 0.452 & {\rm N} & 5.219 & 6.288 & {\rm N} \\
1 & 0.203 & 0.203 & {\rm N} & 0.813 & 0.794 & {\rm N} & 4.576 &  5.420 & & 0.026 & 0.026 & {\rm N} & 0.209 & 0.216 & {\rm N} & 3.176 & 4.354 & {\rm N} \\
3 & 0.219 & 0.220 & {\rm N} & 0.779 & 0.792 & {\rm N} & 4.048 &  4.267 & & 0.026 & 0.025 & {\rm N} & 0.209 & 0.212 & {\rm N} & 3.018 & 4.196 & \\
10 & 3.832 & 2.201 & & 4.974 & 3.642 & & 7.092 &  6.002 & & 3.597 & 0.640 & {\rm P} & 3.620 & 1.665 & {\rm P} & 5.198 & 4.871 & \\

 \hline
 {\rm T} & \multicolumn{3}{c|}{m=2^{2}} & \multicolumn{3}{c|}{m=2^{4}} & \multicolumn{3}{c||}{m=2^{6}} & \multicolumn{3}{c|}{m=2^{3}} & \multicolumn{3}{c|}{m=2^{6}} & \multicolumn{3}{c|}{m=2^{9}}\\ \hline
1/2 & 0.185 & 0.195 & {\rm N} & 0.906 & 2.091 & {\rm N} & 1.770 &  4.144 & {\rm N} & 0.024 & 0.026 & {\rm N} & 0.192 & 0.401 & {\rm N} & 4.522 &  6.006 & \\ 
1 & 0.205 & 0.198 & {\rm N} & 0.834 & 1.064 & {\rm N} & 4.740 &  5.351 & & 0.027 & 0.027 & {\rm N} & 0.207 & 0.210 & {\rm N} & 2.401 & 3.655 & {\rm N} \\ 
3 & 0.199 & 0.182 & {\rm N} & 0.719 & 0.692 & {\rm N} & 3.792 &  4.058 & & 0.025 & 0.025 & {\rm N} & 0.209 & 0.211 & {\rm N} & 1.677 & 1.729 & {\rm N} \\ 
10 & 0.214 & 0.195 & {\rm N} & 0.721 & 0.692 & {\rm N} & 4.010 &  3.917 & & 0.026 & 0.025 & {\rm N} & 1.183 & 1.138 & & 1.849 & 1.879 & \\ 
 
 \hline
 {\rm Burr} & \multicolumn{3}{c|}{m=2^{2}} & \multicolumn{3}{c|}{m=2^{4}} & \multicolumn{3}{c||}{m=2^{6}} & \multicolumn{3}{c|}{m=2^{3}} & \multicolumn{3}{c|}{m=2^{6}} & \multicolumn{3}{c|}{m=2^{9}}\\ \hline
1/2,1/2  & 0.216 & 0.484 & {\rm N} & 0.558 & 1.782 & {\rm N} & 0.951 &  3.930 & {\rm N} & 0.532 & 3.652 & {\rm N} & 1.193 & 1.992 & & 4.535 &  4.637 & \\ 
  1,1/2  & 0.211 & 0.744 & {\rm N} & 1.225 & 2.815 & {\rm N} & 2.606 &  6.231 & {\rm N} & 0.024 & 0.027 & {\rm N} & 0.216 & 0.486 & {\rm N} & 3.747 &  5.308 & \\ 
  3,1/2  & 0.202 & 0.195 & {\rm N} & 0.781 & 0.790 & {\rm N} & 4.211 &  4.595 & & 0.026 & 0.029 & {\rm N} & 0.203 & 0.202 & {\rm N} & 1.775 &  2.077 & {\rm N} \\ 
  1/2,1  & 0.188 & 0.240 & {\rm N} & 1.104 & 2.451 & {\rm N} & 4.668 &  5.842 & & 0.024 & 0.028 & {\rm N} & 0.177 & 0.275 & {\rm N} & 3.037 &  5.230 & \\ 
  1,1  & 0.191 & 0.194 & {\rm N} & 0.751 & 0.841 & {\rm N} & 4.912 &  5.257 & & 0.026 & 0.027 & {\rm N} & 0.206 & 0.220 & {\rm N} & 2.635 &  3.764 & {\rm N} \\ 
  3,1  & 0.197 & 0.196 & {\rm N} & 0.734 & 0.780 & {\rm N} & 3.983 &  4.112 & & 0.025 & 0.024 & {\rm N} & 0.196 & 0.182 & {\rm N} & 1.686 & 1.834 & {\rm N} \\ 
  1/2,3  & 0.201 & 0.199 & {\rm N} & 0.776 & 0.807 & {\rm N} & 4.473 &  4.274 & & 0.027 & 0.028 & {\rm N} & 0.205 & 0.206 & {\rm N} & 1.825 & 1.859 & {\rm N} \\ 
  1,3  & 0.205 & 0.198 & {\rm N} & 0.735 & 0.778 & {\rm N} & 3.786 &  3.788 & & 0.026 & 0.026 & {\rm N} & 0.205 & 0.205 & {\rm N} & 1.578 &  1.590 & {\rm N} \\ 
  3,3  & 0.231 & 0.318 & {\rm N} & 3.348 & 2.570 & & 5.574 & 4.817 & & 2.948 & 0.771 & & 3.464 & 1.602 & {\rm P} & 4.935 & 4.299 & \\ 

 \hline
 {\rm Fr{e}chet} & \multicolumn{3}{c|}{m=2^{2}} & \multicolumn{3}{c|}{m=2^{4}} & \multicolumn{3}{c||}{m=2^{6}} & \multicolumn{3}{c|}{m=2^{3}} & \multicolumn{3}{c|}{m=2^{6}} & \multicolumn{3}{c|}{m=2^{9}}\\ \hline
5 & 0.234 & 1.164 & {\rm N} & 1.146 & 4.302 & {\rm N} & 4.476 & 6.136 & & 1.260 & 5.665 & {\rm N} & 0.371 & 0.571 & {\rm N} & 2.820 & 4.392 & \\ 
  2 & 0.220 & 1.070 & {\rm N} & 1.153 & 2.786 & {\rm N} & 5.429 & 6.843 & & 0.023 & 0.026 & {\rm N} & 0.190 & 0.445 & {\rm N} & 3.072 & 4.688 & \\ 
  1 & 0.208 & 0.215 & {\rm N} & 0.801 & 0.878 & {\rm N} & 4.831 & 5.110 & & 0.026 & 0.028 & {\rm N} & 0.205 & 0.211 & {\rm N} & 2.324 & 3.246 & {\rm N} \\ 
  1/2 & 0.226 & 0.236 & {\rm N} & 0.828 & 0.829 & {\rm N} & 3.970 & 3.980 & & 0.026 & 0.026 & {\rm N} & 0.201 & 0.206 & {\rm N} & 1.741 & 1.851 & {\rm N} \\ 
  1/4 & 0.217 & 0.206 & {\rm N} & 0.802 & 0.781 & {\rm N} & 2.864 & 3.489 & & 0.025 & 0.024 & {\rm N} & 0.195 & 0.196 & {\rm N} & 1.833 & 2.236 & {\rm N} \\ 
  
 \hline
 {\rm Weibull} & \multicolumn{3}{c|}{m=2^{2}} & \multicolumn{3}{c|}{m=2^{4}} & \multicolumn{3}{c||}{m=2^{6}} & \multicolumn{3}{c|}{m=2^{3}} & \multicolumn{3}{c|}{m=2^{6}} & \multicolumn{3}{c|}{m=2^{9}}\\ \hline
1/2 & 0.561 & 0.391 & & 1.128 & 1.054 & & 5.826 & 3.503 & {\rm P} & 0.479 & 0.144 & & 4.023 & 0.803 & {\rm P} & 15.54 & 1.703 & {\rm P} \\ 
1/1 & 0.282 & 0.358 & & 1.761 & 1.511 & {\rm P} & 7.516 &  4.369 & {\rm P} & 0.753 & 0.595 & & 6.943 & 1.108 & {\rm P} & 17.61 & 1.638 & {\rm P} \\ 
3 & 2.121 & 0.924 & {\rm P} & 2.332 & 1.775 & {\rm P} & 8.023 &  4.376 & {\rm P} & 2.073 & 0.298 & {\rm P} & 8.291 & 1.202 & {\rm P} & 18.02 & 2.272 & {\rm P} \\ 
10 & 2.717 & 1.028 & {\rm P} & 2.529 & 1.868 & {\rm P} & 8.651 &  4.760 & {\rm P} & 2.560 & 0.356 & {\rm P} & 8.757 & 1.237 & {\rm P} & 16.79 & 4.632 & {\rm P} \\ 
 
 \hline
 {\rm inv. Burr} & \multicolumn{3}{c|}{m=2^{2}} & \multicolumn{3}{c|}{m=2^{4}} & \multicolumn{3}{c||}{m=2^{6}} & \multicolumn{3}{c|}{m=2^{3}} & \multicolumn{3}{c|}{m=2^{6}} & \multicolumn{3}{c|}{m=2^{9}}\\ \hline
  3,2  & 0.196 & 0.202 & {\rm N} & 1.358 & 1.540 & & 4.766 &  5.020 & & 0.118 & 0.273 & {\rm N} & 1.729 & 1.124 & {\rm P} & 4.858 &  5.125 & {\rm P} \\ 
  1,2  & 0.182 & 0.205 & {\rm N} & 1.813 & 2.783 & {\rm P} & 6.097 &  6.062 & {\rm P}  & 0.678 & 0.287 & & 0.999 & 0.836 & {\rm P} & 5.056 &  7.057 & {\rm P} \\
   1/2,2  & 1.412 & 4.658 & & 4.928 & 8.703 & {\rm P} & 7.347 &  10.50 & {\rm P} & 0.984 & 4.283 & {\rm P} & 2.871 & 3.729 & {\rm P} & 5.356 &  6.533 & {\rm P} \\ 
 3,1  & 0.500 & 2.364 & {\rm N} & 0.806 & 1.287 & & 3.245 &  4.743 & & 0.252 & 1.584 & {\rm N} & 1.032 & 0.865 & {\rm P} & 4.592 &  5.931 & {\rm P} \\ 
 1,1  & 1.363 & 4.243 & & 2.329 & 4.328 & {\rm P} & 6.386 &  7.855 & {\rm P} & 0.259 & 0.744 & & 3.510 & 4.222 & {\rm P} & 6.516 &  9.926 & {\rm P} \\ 
 1/2,1  & 3.646 & 7.507 & & 7.646 & 11.11 & {\rm P} & 8.390 &  12.32 & {\rm P} & 5.408 & 10.11 & {\rm P} & 4.863 & 3.874 & {\rm P} & 6.949 &  7.682 & {\rm P} \\ 
3,1/3  & 1.273 & 5.067 & {\rm N} & 0.525 & 1.734 & {\rm N} & 1.611 &  3.971 & {\rm N} & 1.412 & 4.658 & & 3.772 & 7.368 & {\rm P} & 5.697 &  8.013 & {\rm P} \\ 
  1,1/3  & 1.537 & 5.052 & {\rm N} & 2.950 & 5.307 & {\rm P} & 6.851 &  9.248 & {\rm P} & 0.138 & 0.575 & {\rm N} & 5.694 & 4.964 & {\rm P} & 8.365 & 11.31 & {\rm P} \\ 
  1/2,1/3  & 4.817 & 8.808 & & 9.287 & 12.40 & {\rm N} & 8.443 & 13.98 & & 6.113 & 11.41 & {\rm P} & 7.026 & 3.756 & {\rm P} & 8.663 & 8.683 & {\rm P} \\ 
  
 \hline
 \text{non-MDA} & \multicolumn{3}{c|}{m=2^{2}} & \multicolumn{3}{c|}{m=2^{4}} & \multicolumn{3}{c||}{m=2^{6}} & \multicolumn{3}{c|}{m=2^{3}} & \multicolumn{3}{c|}{m=2^{6}} & \multicolumn{3}{c|}{m=2^{9}}\\ \hline
 & 0.207 & 0.216 & {\rm N} & 1.090 & 1.258 & {\rm N} & 2.155 & 3.550 & {\rm N} & 0.031 & 0.036 & & 1.218 & 0.732 & & 1.789 & 2.083 & \\

\hline
 \end{tabu}$$
}
\end{table}
\end{landscape}

\begin{landscape}
\begin{table}
\caption{{Scaled MISE values ($\times 100$) and sd values ($\times 100$) of the proposed approach with AD metric}}
{\fontsize{6pt}{6pt}\selectfont
$$\begin{tabu}[c]{|c||ccc|ccc|ccc||ccc|ccc|ccc|}
 \hline
 & \multicolumn{9}{c||}{n=2^8} & \multicolumn{9}{c|}{n=2^{12}} \\ \hline
 & {\rm MISE} & {\rm sd} & {\rm P/N} & {\rm MISE} & {\rm sd} & {\rm P/N} & {\rm MISE} & {\rm sd} & {\rm P/N} & {\rm MISE} & {\rm sd} & {\rm P/N} & {\rm MISE} & {\rm sd} & {\rm P/N} & {\rm MISE} & {\rm sd} & {\rm P/N} \\ \hline\hline 
 {\rm Pareto} & \multicolumn{3}{c|}{m=2^{2}} & \multicolumn{3}{c|}{m=2^{4}} & \multicolumn{3}{c||}{m=2^{6}} & \multicolumn{3}{c|}{m=2^{3}} & \multicolumn{3}{c|}{m=2^{6}} & \multicolumn{3}{c|}{m=2^{9}} \\ \hline
1/2 & 0.163 & 0.256 &  & 0.760 & 1.666 &  & 1.584 &  4.063 & {\rm N} & 0.041 & 0.184 &  & 0.176 & 0.296 &  & 3.742 &  5.208 &  \\ 
1 & 0.169 & 0.206 & {\rm P} & 0.758 & 1.119 & {\rm P} & 4.143 &  4.601 &  & 0.028 & 0.044 & {\rm P} & 0.174 & 0.212 & {\rm P} & 2.252 &  3.021 &  \\ 
3 & 0.195 & 0.228 & {\rm P} & 0.830 & 1.032 & {\rm P} & 3.988 &  4.131 &  & 0.025 & 0.028 & {\rm P} & 0.181 & 0.206 & {\rm P} & 1.792 &  1.980 & {\rm P} \\ 
10 & 0.210 & 0.269 & {\rm P} & 0.862 & 1.079 & {\rm P} & 6.281 &  5.161 & {\rm P} & 0.026 & 0.029 & {\rm P} & 0.201 & 0.239 & {\rm P} & 1.986 &  2.280 & {\rm P} \\ 

 \hline
 {\rm T} & \multicolumn{3}{c|}{m=2^{2}} & \multicolumn{3}{c|}{m=2^{4}} & \multicolumn{3}{c||}{m=2^{6}} & \multicolumn{3}{c|}{m=2^{3}} & \multicolumn{3}{c|}{m=2^{6}} & \multicolumn{3}{c|}{m=2^{9}}\\ \hline
1/2 & 0.325 & 0.869 & {\rm P} & 0.827 & 2.043 &  & 1.534 &  4.059 & {\rm N} & 0.162 & 0.490 & {\rm P} & 0.168 & 0.214 &  & 3.161 &  4.125 &  \\ 
1 & 0.261 & 0.333 & {\rm P} & 0.769 & 1.322 & {\rm P} & 5.954 &  5.960 &  & 0.044 & 0.131 & {\rm P} & 0.174 & 0.295 & {\rm P} & 1.959 &  2.754 &  \\ 
3 & 0.241 & 0.261 & {\rm P} & 0.820 & 1.003 & {\rm P} & 4.504 &  4.712 &  & 0.045 & 0.045 & {\rm P} & 0.199 & 0.238 & {\rm P} & 1.953 &  2.170 & {\rm P} \\ 
10 & 0.214 & 0.240 & {\rm P} & 0.994 & 1.161 & {\rm P} & 4.702 &  4.401 &  & 0.037 & 0.035 & {\rm P} & 0.197 & 0.237 & {\rm P} & 1.991 &  2.114 & {\rm P} \\ 

 \hline
 {\rm Burr} & \multicolumn{3}{c|}{m=2^{2}} & \multicolumn{3}{c|}{m=2^{4}} & \multicolumn{3}{c||}{m=2^{6}} & \multicolumn{3}{c|}{m=2^{3}} & \multicolumn{3}{c|}{m=2^{6}} & \multicolumn{3}{c|}{m=2^{9}}\\ \hline
1/2,1/2  & 0.144 & 0.188 & {\rm N} & 0.422 & 0.943 & {\rm N} & 0.837 &  2.715 & {\rm N} & 0.033 & 0.110 & {\rm N} & 0.326 & 0.547 &  & 4.551 &  4.783 &  \\ 
  1,1/2  & 0.193 & 0.378 &  & 0.787 & 1.541 &  & 1.873 &  3.726 & {\rm N} & 0.037 & 0.139 &  & 0.176 & 0.224 &  & 3.232 &  4.685 &  \\ 
  3,1/2  & 0.180 & 0.222 & {\rm P} & 0.722 & 0.846 & {\rm P} & 5.010 &  4.890 &  & 0.029 & 0.035 & {\rm P} & 0.198 & 0.246 & {\rm P} & 1.759 &  2.206 &  \\ 
  1/2,1  & 0.183 & 0.329 &  & 0.768 & 1.441 &  & 4.749 &  4.730 &  & 0.064 & 0.316 &  & 0.169 & 0.234 &  & 2.046 &  3.934 &  \\ 
  1,1  & 0.181 & 0.224 & {\rm P} & 0.858 & 1.511 & {\rm P} & 4.492 &  4.214 &  & 0.031 & 0.067 & {\rm P} & 0.191 & 0.300 & {\rm P} & 1.840 &  2.323 &  \\ 
  3,1  & 0.201 & 0.240 & {\rm P} & 0.800 & 0.926 & {\rm P} & 4.381 &  3.984 &  & 0.030 & 0.036 & {\rm P} & 0.190 & 0.221 & {\rm P} & 1.856 &  2.075 & {\rm P} \\ 
  1/2,3  & 0.172 & 0.202 & {\rm P} & 0.793 & 1.294 & {\rm P} & 5.667 &  4.767 & {\rm P} & 0.027 & 0.025 & {\rm P} & 0.189 & 0.265 & {\rm P} & 1.712 &  1.983 & {\rm P} \\ 
  1,3  & 0.197 & 0.245 & {\rm P} & 0.837 & 1.004 & {\rm P} & 4.535 &  4.134 &  & 0.026 & 0.031 & {\rm P} & 0.196 & 0.233 & {\rm P} & 1.961 &  2.352 & {\rm P} \\ 
  3,3  & 0.215 & 0.273 & {\rm P} & 0.895 & 0.997 & {\rm P} & 6.186 &  4.887 & {\rm P} & 0.034 & 0.037 & {\rm P} & 0.207 & 0.240 & {\rm P} & 1.967 &  2.264 & {\rm P} \\ 

 \hline
 {\rm Fr{e}chet} & \multicolumn{3}{c|}{m=2^{2}} & \multicolumn{3}{c|}{m=2^{4}} & \multicolumn{3}{c||}{m=2^{6}} & \multicolumn{3}{c|}{m=2^{3}} & \multicolumn{3}{c|}{m=2^{6}} & \multicolumn{3}{c|}{m=2^{9}} \\ \hline
5 & 0.146 & 0.188 & {\rm N} & 0.561 & 0.843 & {\rm N} & 3.795 &  4.598 &  & 0.025 & 0.064 & {\rm N} & 0.153 & 0.188 & {\rm N} & 1.739 &  2.604 &  \\
2 & 0.230 & 0.496 & {\rm P} & 0.861 & 1.987 &  & 4.698 &  5.044 &  & 0.092 & 0.244 &  & 0.184 & 0.292 &  & 1.844 &  2.658 &  \\ 
1 & 0.181 & 0.233 & {\rm P} & 0.756 & 1.301 & {\rm P} & 4.614 &  4.788 &  & 0.029 & 0.060 & {\rm P} & 0.198 & 0.400 & {\rm P} & 1.853 &  2.278 &  \\ 
  1/2  & 0.185 & 0.211 & {\rm P} & 0.804 & 0.915 & {\rm P} & 4.761 &  4.544 &  & 0.029 & 0.037 & {\rm P} & 0.205 & 0.237 & {\rm P} & 1.790 &  2.074 & {\rm P} \\ 
  1/4 & 0.211 & 0.271 & {\rm P} & 0.895 & 1.053 & {\rm P} & 2.985 &  3.761 &  & 0.029 & 0.035 & {\rm P} & 0.196 & 0.225 & {\rm P} & 2.059 &  2.523 & {\rm P} \\ 
 
 \hline
 {\rm Weibull} & \multicolumn{3}{c|}{m=2^{2}} & \multicolumn{3}{c|}{m=2^{4}} & \multicolumn{3}{c||}{m=2^{6}} & \multicolumn{3}{c|}{m=2^{3}} & \multicolumn{3}{c|}{m=2^{6}} & \multicolumn{3}{c|}{m=2^{9}}\\ \hline
1/2 & 2.585 & 0.636 & {\rm P} & 3.016 & 1.619 & {\rm P} & 8.464 &  3.387 & {\rm P} & 2.641 & 0.272 & {\rm P} & 7.688 & 0.854 & {\rm P} & 16.30 &  1.179 & {\rm P} \\ 
1 & 3.268 & 0.976 & {\rm P} & 3.343 & 1.934 & {\rm P} & 9.175 &  3.937 & {\rm P} & 3.046 & 0.340 & {\rm P} & 8.773 & 0.987 & {\rm P} & 18.00 &  1.504 & {\rm P} \\ 
3 & 3.443 & 1.034 & {\rm P} & 3.581 & 2.118 & {\rm P} & 9.344 &  4.342 & {\rm P} & 3.324 & 0.353 & {\rm P} & 9.424 & 1.103 & {\rm P} & 18.22 &  2.249 & {\rm P} \\ 
10 & 3.666 & 1.045 & {\rm P} & 3.713 & 2.160 & {\rm P} & 9.499 &  4.379 & {\rm P} & 3.513 & 0.367 & {\rm P} & 9.607 & 1.109 & {\rm P} & 17.01 &  4.395 & {\rm P} \\ 
 
\hline
 {\rm inv. Burr} & \multicolumn{3}{c|}{m=2^{2}} & \multicolumn{3}{c|}{m=2^{4}} & \multicolumn{3}{c||}{m=2^{6}} & \multicolumn{3}{c|}{m=2^{3}} & \multicolumn{3}{c|}{m=2^{6}} & \multicolumn{3}{c|}{m=2^{9}}\\ \hline
  3,2  & 0.245 & 0.258 & {\rm P} & 0.887 & 1.010 & {\rm P} & 5.143 &  4.559 &  & 0.066 & 0.094 & {\rm P} & 0.208 & 0.248 & {\rm P} & 2.141 &  2.431 & {\rm P} \\ 
    1,2  & 0.225 & 0.239 & {\rm P} & 0.915 & 1.047 & {\rm P} & 5.488 &  5.362 & {\rm P} & 0.040 & 0.093 & {\rm P} & 0.212 & 0.253 & {\rm P} & 2.427 &  3.188 & {\rm P} \\ 
    1/2,2  & 1.008 & 3.091 & {\rm P} & 2.992 & 4.597 & {\rm P} & 5.488 &  6.411 & {\rm P} & 0.401 & 1.524 & {\rm P} & 2.425 & 3.535 & {\rm P} & 3.982 &  5.513 & {\rm P} \\ 
  3,1  & 0.321 & 1.317 & {\rm P} & 0.776 & 1.035 & {\rm P} & 3.320 &  4.109 &  & 0.420 & 1.649 & {\rm P} & 0.212 & 0.247 & {\rm P} & 2.581 &  2.825 & {\rm P} \\   
    1,1  & 0.635 & 1.593 & {\rm P} & 1.371 & 2.333 & {\rm P} & 5.409 &  6.318 & {\rm P} & 0.183 & 0.474 & {\rm P} & 2.868 & 3.845 & {\rm P} & 2.858 &  3.731 & {\rm P} \\   
      1/2,1  & 2.722 & 5.041 & {\rm P} & 4.445 & 5.415 & {\rm P} & 5.528 &  6.783 & {\rm P} & 1.715 & 3.255 & {\rm P} & 4.461 & 4.060 & {\rm P} & 5.151 &  6.770 &  \\ 
3,1/3  & 0.084 & 1.255 & {\rm N} & 0.244 & 0.735 & {\rm N} & 1.442 &  2.977 & {\rm N} & 1.008 & 3.091 & {\rm P} & 2.330 & 3.726 & {\rm P} & 4.331 &  5.362 & {\rm P} \\   
  1,1/3  & 0.752 & 1.900 & {\rm P} & 1.416 & 2.941 & {\rm P} & 5.567 &  7.227 &  & 0.088 & 0.516 & {\rm P} & 5.035 & 4.038 & {\rm P} & 3.298 &  4.531 &  \\ 
  1/2,1/3  & 3.658 & 5.804 & {\rm P} & 5.759 & 6.538 & {\rm P} & 6.414 &  8.202 &  & 0.991 & 2.895 & {\rm P} & 7.482 & 4.087 & {\rm P} & 9.081 &  8.300 &  \\ 

\hline
 
 \text{non-MDA} & \multicolumn{3}{c|}{m=2^{2}} & \multicolumn{3}{c|}{m=2^{4}} & \multicolumn{3}{c||}{m=2^{6}} & \multicolumn{3}{c|}{m=2^{3}} & \multicolumn{3}{c|}{m=2^{6}} & \multicolumn{3}{c|}{m=2^{9}}\\ \hline
 
& 0.336 & 0.176 & {\rm P} & 2.169 & 1.598 & {\rm P} & 2.099 & 3.247 & {\rm P} & 0.610 & 0.089 & {\rm P} & 0.839 &  0.493 & {\rm P} & 1.925 & 2.152 & {\rm P} \\

\hline
 \end{tabu}$$
}
\end{table}
\end{landscape}

\begin{table}[h]
\caption{{Estimated mixing parameter $p$ and sd values of the naive mixture approach with ALbw function}}
{\fontsize{6pt}{6pt}\selectfont
$$\begin{tabu}[c]{|c||cc|cc|cc||cc|cc|cc|}
 \hline
 & \multicolumn{6}{c||}{n=2^8} & \multicolumn{6}{c|}{n=2^{12}} \\ \hline
 & p & {\rm sd} & p & {\rm sd} & p & {\rm sd} & p & {\rm sd} & p & {\rm sd} & p & {\rm sd}\\ \hline\hline 
 {\rm Pareto} & \multicolumn{2}{c|}{m=2^{2}} & \multicolumn{2}{c|}{m=2^{4}} & \multicolumn{2}{c||}{m=2^{6}} & \multicolumn{2}{c|}{m=2^{3}} & \multicolumn{2}{c|}{m=2^{6}} & \multicolumn{2}{c|}{m=2^{9}} \\ \hline
1/2 & 0.987 & 0.042 & 0.949 & 0.201 & 0.788 & 0.391 & 0.991 & 0.033 & 0.959 & 0.149 & 0.795 & 0.393 \\ 
1 & 0.688 & 0.258 & 0.665 & 0.459 & 0.688 & 0.463 & 0.360 & 0.235 & 0.605 & 0.485 & 0.689 & 0.463 \\ 
3 & 0.134 & 0.167 & 0.326 & 0.468 & 0.541 & 0.498 & 0.002 & 0.005 & 0.215 & 0.411 & 0.509 & 0.500 \\ 
10 & 0.040 & 0.069 & 0.179 & 0.382 & 0.482 & 0.498 & 0.000 & 0.000 & 0.056 & 0.230 & 0.400 & 0.490 \\ 

 \hline
 {\rm T} & \multicolumn{2}{c|}{m=2^{2}} & \multicolumn{2}{c|}{m=2^{4}} & \multicolumn{2}{c||}{m=2^{6}} & \multicolumn{2}{c|}{m=2^{3}} & \multicolumn{2}{c|}{m=2^{6}} & \multicolumn{2}{c|}{m=2^{9}}\\ \hline
1/2 & 1.000 & 0.000 & 0.999 & 0.032 & 0.896 & 0.278 & 1.000 & 0.000 & 1.000 & 0.000 & 0.895 & 0.285 \\ 
1 & 0.998 & 0.038 & 0.944 & 0.226 & 0.695 & 0.456 & 1.000 & 0.000 & 0.994 & 0.077 & 0.747 & 0.434 \\ 
3 & 0.350 & 0.440 & 0.531 & 0.498 & 0.528 & 0.498 & 0.985 & 0.121 & 0.842 & 0.365 & 0.542 & 0.498 \\ 
10 & 0.022 & 0.070 & 0.205 & 0.403 & 0.463 & 0.497 & 0.000 & 0.001 & 0.169 & 0.375 & 0.407 & 0.491 \\ 
 
 \hline
 {\rm Burr} & \multicolumn{2}{c|}{m=2^{2}} & \multicolumn{2}{c|}{m=2^{4}} & \multicolumn{2}{c||}{m=2^{6}} & \multicolumn{2}{c|}{m=2^{3}} & \multicolumn{2}{c|}{m=2^{6}} & \multicolumn{2}{c|}{m=2^{9}}\\ \hline
1/2,1/2 & 0.997 & 0.008 & 0.998 & 0.018 & 0.904 & 0.253 & 0.948 & 0.213 & 0.767 & 0.423 & 0.260 & 0.431 \\ 
 1,1/2 & 0.988 & 0.037 & 0.962 & 0.169 & 0.815 & 0.368 & 0.993 & 0.028 & 0.955 & 0.148 & 0.547 & 0.486 \\ 
 3,1/2 & 0.556 & 0.391 & 0.561 & 0.496 & 0.615 & 0.485 & 0.620 & 0.457 & 0.520 & 0.500 & 0.612 & 0.487 \\ 
 1/2,1 & 0.988 & 0.034 & 0.964 & 0.155 & 0.803 & 0.379 & 0.990 & 0.063 & 0.961 & 0.138 & 0.670 & 0.461 \\ 
 1,1 & 0.698 & 0.262 & 0.636 & 0.468 & 0.698 & 0.459 & 0.354 & 0.233 & 0.604 & 0.485 & 0.666 & 0.470 \\ 
 3,1 & 0.161 & 0.275 & 0.404 & 0.490 & 0.540 & 0.497 & 0.036 & 0.164 & 0.463 & 0.499 & 0.585 & 0.493 \\ 
 1/2,3 & 0.790 & 0.191 & 0.445 & 0.481 & 0.617 & 0.486 & 0.409 & 0.136 & 0.301 & 0.458 & 0.620 & 0.486 \\ 
 1,3 & 0.133 & 0.143 & 0.314 & 0.464 & 0.535 & 0.498 & 0.001 & 0.005 & 0.251 & 0.434 & 0.550 & 0.498 \\ 
 3,3 & 0.024 & 0.049 & 0.221 & 0.414 & 0.455 & 0.496 & 0.000 & 0.002 & 0.163 & 0.369 & 0.476 & 0.500 \\ 

 \hline
 {\rm Fr{e}chet} & \multicolumn{2}{c|}{m=2^{2}} & \multicolumn{2}{c|}{m=2^{4}} & \multicolumn{2}{c||}{m=2^{6}} & \multicolumn{2}{c|}{m=2^{3}} & \multicolumn{2}{c|}{m=2^{6}} & \multicolumn{2}{c|}{m=2^{9}}\\ \hline
5 & 0.994 & 0.01 & 0.998 & 0.01 & 0.938 & 0.202 & 0.999 & 0.001 & 1.000 & 0.001 & 0.973 & 0.134 \\ 
2 & 0.988 & 0.039 & 0.953 & 0.188 & 0.833 & 0.356 & 0.992 & 0.028 & 0.961 & 0.139 & 0.813 & 0.378 \\ 
1 & 0.730 & 0.272 & 0.629 & 0.473 & 0.667 & 0.470 & 0.537 & 0.355 & 0.562 & 0.493 & 0.689 & 0.463 \\ 
 1/2 & 0.261 & 0.297 & 0.451 & 0.497 & 0.595 & 0.490 & 0.052 & 0.168 & 0.436 & 0.496 & 0.564 & 0.496 \\ 
1/4 & 0.100 & 0.195 & 0.333 & 0.471 & 0.506 & 0.498 & 0.003 & 0.006 & 0.267 & 0.443 & 0.476 & 0.500 \\ 
 
 \hline
 {\rm Weibull} & \multicolumn{2}{c|}{m=2^{2}} & \multicolumn{2}{c|}{m=2^{4}} & \multicolumn{2}{c||}{m=2^{6}} & \multicolumn{2}{c|}{m=2^{3}} & \multicolumn{2}{c|}{m=2^{6}} & \multicolumn{2}{c|}{m=2^{9}}\\ \hline
1/2 & 0.434 & 0.219 & 0.180 & 0.384 & 0.491 & 0.499 & 0.004 & 0.012 & 0.020 & 0.140 & 0.400 & 0.490 \\ 
1 & 0.023 & 0.047 & 0.138 & 0.345 & 0.434 & 0.494 & 0.000 & 0.000 & 0.013 & 0.113 & 0.341 & 0.474 \\ 
3 & 0.011 & 0.025 & 0.088 & 0.282 & 0.392 & 0.486 & 0.000 & 0.001 & 0.008 & 0.089 & 0.331 & 0.471 \\ 
10 & 0.009 & 0.021 & 0.064 & 0.241 & 0.367 & 0.480 & 0.000 & 0.001 & 0.010 & 0.100 & 0.318 & 0.466 \\ 
 
 \hline
 {\rm inv. Burr} & \multicolumn{2}{c|}{m=2^{2}} & \multicolumn{2}{c|}{m=2^{4}} & \multicolumn{2}{c||}{m=2^{6}} & \multicolumn{2}{c|}{m=2^{3}} & \multicolumn{2}{c|}{m=2^{6}} & \multicolumn{2}{c|}{m=2^{9}}\\ \hline
 3,2 & 0.021 & 0.055 & 0.043 & 0.184 & 0.375 & 0.473 & 0.000 & 0.001 & 0.000 & 0.003 & 0.261 & 0.439 \\ 
1,2 & 0.019 & 0.055 & 0.018 & 0.107 & 0.302 & 0.446 & 0.003 & 0.046 & 0.001 & 0.011 & 0.089 & 0.282 \\ 
 1/2,2 & 0.034 & 0.050 & 0.015 & 0.043 & 0.256 & 0.399 & 0.044 & 0.203 & 0.001 & 0.017 & 0.046 & 0.162 \\  
 3,1 & 0.451 & 0.221 & 0.322 & 0.208 & 0.642 & 0.364 & 0.316 & 0.232 & 0.197 & 0.096 & 0.513 & 0.392 \\ 
 1,1 & 0.496 & 0.223 & 0.317 & 0.195 & 0.457 & 0.379 & 0.411 & 0.298 & 0.285 & 0.134 & 0.392 & 0.296 \\
 1/2,1 & 0.536 & 0.191 & 0.311 & 0.177 & 0.249 & 0.321 & 0.371 & 0.214 & 0.322 & 0.285 & 0.109 & 0.180 \\  
3,1/3 & 0.853 & 0.140 & 0.518 & 0.198 & 0.360 & 0.340 & 0.424 & 0.101 & 0.324 & 0.177 & 0.228 & 0.205 \\ 
 1,1/3 & 0.891 & 0.145 & 0.364 & 0.164 & 0.154 & 0.241 & 0.428 & 0.058 & 0.659 & 0.402 & 0.043 & 0.093 \\ 
 1/2,1/3 & 0.790 & 0.184 & 0.385 & 0.279 & 0.076 & 0.185 & 0.177 & 0.246 & 0.779 & 0.257 & 0.043 & 0.173 \\

 \hline
 
 \text{non-MDA} & \multicolumn{2}{c|}{m=2^{2}} & \multicolumn{2}{c|}{m=2^{4}} & \multicolumn{2}{c||}{m=2^{6}} & \multicolumn{2}{c|}{m=2^{3}} & \multicolumn{2}{c|}{m=2^{6}} & \multicolumn{2}{c|}{m=2^{9}}\\ \hline
 & 0.206 & 0.264 & 0.263 & 0.440 & 0.300 & 0.458 & 0.016 & 0.122 & 0.001 & 0.032 & 0.548 & 0.498 \\

\hline
 \end{tabu}$$
}
\end{table}

\clearpage

\begin{table}[h]
\caption{{Estimated mixing parameter $p$ and sd values of the naive mixture approach with CVbw function}}
{\fontsize{6pt}{6pt}\selectfont
$$\begin{tabu}[c]{|c||cc|cc|cc||cc|cc|cc|}
 \hline
 & \multicolumn{6}{c||}{n=2^8} & \multicolumn{6}{c|}{n=2^{12}} \\ \hline
 & p & {\rm sd} & p & {\rm sd} & p & {\rm sd} & p & {\rm sd} & p & {\rm sd} & p & {\rm sd}\\ \hline\hline 
 {\rm Pareto} & \multicolumn{2}{c|}{m=2^{2}} & \multicolumn{2}{c|}{m=2^{4}} & \multicolumn{2}{c||}{m=2^{6}} & \multicolumn{2}{c|}{m=2^{3}} & \multicolumn{2}{c|}{m=2^{6}} & \multicolumn{2}{c|}{m=2^{9}} \\ \hline
1/2 & 0.982 & 0.057 & 0.953 & 0.189 & 0.795 & 0.391 & 0.994 & 0.021 & 0.970 & 0.150 & 0.503 & 0.492 \\
  1 & 0.658 & 0.275 & 0.626 & 0.475 & 0.703 & 0.457 & 0.676 & 0.210 & 0.661 & 0.428 & 0.540 & 0.501 \\
  3 & 0.053 & 0.118 & 0.326 & 0.469 & 0.517 & 0.499 & 0.037 & 0.020 & 0.130 & 0.338 & 0.500 & 0.502 \\ 
  10 & 0.026 & 0.056 & 0.193 & 0.395 & 0.475 & 0.497 & 0.000 & 0.002 & 0.055 & 0.227 & 0.395 & 0.481 \\

 \hline
 {\rm T} & \multicolumn{2}{c|}{m=2^{2}} & \multicolumn{2}{c|}{m=2^{4}} & \multicolumn{2}{c||}{m=2^{6}} & \multicolumn{2}{c|}{m=2^{3}} & \multicolumn{2}{c|}{m=2^{6}} & \multicolumn{2}{c|}{m=2^{9}}\\ \hline
1/2 & 1.000 & 0.000 & 0.996 & 0.047 & 0.881 & 0.303 & 1.000 & 0.000 & 1.000 & 0.000 & 0.536 & 0.494 \\
  1 & 0.999 & 0.020 & 0.946 & 0.222 & 0.693 & 0.457 & 1.000 & 0.000 & 1.000 & 0.000 & 0.663 & 0.456 \\
  3 & 0.345 & 0.444 & 0.560 & 0.496 & 0.532 & 0.499 & 0.980 & 0.140 & 0.850 & 0.359 & 0.570 & 0.498 \\
  10 & 0.023 & 0.082 & 0.206 & 0.404 & 0.455 & 0.496 & 0.000 & 0.001 & 0.160 & 0.368 & 0.415 & 0.495 \\
 
 \hline
 {\rm Burr} & \multicolumn{2}{c|}{m=2^{2}} & \multicolumn{2}{c|}{m=2^{4}} & \multicolumn{2}{c||}{m=2^{6}} & \multicolumn{2}{c|}{m=2^{3}} & \multicolumn{2}{c|}{m=2^{6}} & \multicolumn{2}{c|}{m=2^{9}}\\ \hline
1/2,1/2 & 0.873 & 0.311 & 0.670 & 0.468 & 0.216 & 0.371 & 0.937 & 0.238 & 0.700 & 0.461 & 0.245 & 0.426 \\
  1,1/2 & 0.981 & 0.056 & 0.904 & 0.278 & 0.271 & 0.405 & 0.996 & 0.017 & 0.992 & 0.050 & 0.599 & 0.457 \\
  3,1/2 & 0.507 & 0.418 & 0.595 & 0.491 & 0.661 & 0.465 & 0.761 & 0.206 & 0.237 & 0.402 & 0.576 & 0.492 \\
  1/2,1 & 0.975 & 0.067 & 0.891 & 0.291 & 0.571 & 0.477 & 0.985 & 0.101 & 0.981 & 0.111 & 0.453 & 0.484 \\ 
  1,1 & 0.650 & 0.280 & 0.646 & 0.468 & 0.765 & 0.421 & 0.677 & 0.203 & 0.689 & 0.433 & 0.500 & 0.502 \\ 
  3,1 & 0.152 & 0.285 & 0.444 & 0.496 & 0.698 & 0.458 & 0.037 & 0.023 & 0.240 & 0.429 & 0.480 & 0.502 \\ 
  1/2,3 & 0.741 & 0.216 & 0.442 & 0.484 & 0.740 & 0.439 & 0.669 & 0.410 & 0.486 & 0.497 & 0.615 & 0.489 \\
  1,3 & 0.054 & 0.133 & 0.313 & 0.464 & 0.723 & 0.448 & 0.050 & 0.169 & 0.510 & 0.502 & 0.530 & 0.502 \\ 
  3,3 & 0.028 & 0.074 & 0.228 & 0.419 & 0.698 & 0.458 & 0.000 & 0.001 & 0.200 & 0.402 & 0.360 & 0.482 \\ 

 \hline
 {\rm Fr{e}chet} & \multicolumn{2}{c|}{m=2^{2}} & \multicolumn{2}{c|}{m=2^{4}} & \multicolumn{2}{c||}{m=2^{6}} & \multicolumn{2}{c|}{m=2^{3}} & \multicolumn{2}{c|}{m=2^{6}} & \multicolumn{2}{c|}{m=2^{9}}\\ \hline
   5 & 0.979 & 0.021 & 0.959 & 0.178 & 0.669 & 0.436 & 0.867 & 0.347 & 0.705 & 0.456 & 0.135 & 0.343 \\
     2 & 0.981 & 0.050 & 0.891 & 0.291 & 0.671 & 0.462 & 0.999 & 0.006 & 0.952 & 0.188 & 0.402 & 0.481 \\ 
       1 & 0.699 & 0.295 & 0.623 & 0.477 & 0.793 & 0.405 & 0.748 & 0.212 & 0.642 & 0.429 & 0.495 & 0.499 \\ 
 1/2 & 0.228 & 0.338 & 0.472 & 0.499 & 0.711 & 0.453 & 0.130 & 0.204 & 0.420 & 0.496 & 0.600 & 0.492 \\
1/4 & 0.080 & 0.185 & 0.348 & 0.476 & 0.521 & 0.500 & 0.010 & 0.014 & 0.280 & 0.451 & 0.510 & 0.502 \\
 
 \hline
 {\rm Weibull} & \multicolumn{2}{c|}{m=2^{2}} & \multicolumn{2}{c|}{m=2^{4}} & \multicolumn{2}{c||}{m=2^{6}} & \multicolumn{2}{c|}{m=2^{3}} & \multicolumn{2}{c|}{m=2^{6}} & \multicolumn{2}{c|}{m=2^{9}}\\ \hline
1/2 & 0.619 & 0.128 & 0.107 & 0.228 & 0.570 & 0.336 & 0.407 & 0.090 & 0.001 & 0.008 & 0.944 & 0.210 \\
  1 & 0.296 & 0.169 & 0.044 & 0.166 & 0.442 & 0.309 & 0.000 & 0.000 & 0.020 & 0.141 & 0.320 & 0.469 \\
  3 & 0.214 & 0.107 & 0.026 & 0.094 & 0.397 & 0.299 & 0.036 & 0.019 & 0.000 & 0.000 & 0.634 & 0.286 \\ 
  10 & 0.216 & 0.094 & 0.027 & 0.089 & 0.382 & 0.295 & 0.037 & 0.016 & 0.001 & 0.010 & 0.526 &  0.288 \\ 
 
 \hline
 {\rm inv. Burr} & \multicolumn{2}{c|}{m=2^{2}} & \multicolumn{2}{c|}{m=2^{4}} & \multicolumn{2}{c||}{m=2^{6}} & \multicolumn{2}{c|}{m=2^{3}} & \multicolumn{2}{c|}{m=2^{6}} & \multicolumn{2}{c|}{m=2^{9}}\\ \hline
   3,2 & 0.018 & 0.042 & 0.052 & 0.198 & 0.556 & 0.479 & 0.051 & 0.111 & 0.062 & 0.078 & 0.378 & 0.454 \\
  1,2 & 0.015 & 0.039 & 0.053 & 0.164 & 0.504 & 0.476 & 0.050 & 0.094 & 0.063 & 0.096 & 0.106 & 0.236 \\ 
  1/2,2 & 0.133 & 0.264 & 0.211 & 0.277 & 0.313 & 0.424 & 0.216 & 0.389 & 0.356 & 0.391 & 0.116 & 0.235 \\     
  3,1 & 0.305 & 0.259 & 0.217 & 0.218 & 0.437 & 0.455 & 0.472 & 0.181 & 0.291 & 0.093 & 0.406 & 0.423 \\ 
  1,1 & 0.294 & 0.338 & 0.204 & 0.233 & 0.450 & 0.445 & 0.174 & 0.266 & 0.167 & 0.204 & 0.158 & 0.303 \\ 
  1/2,1 & 0.378 & 0.370 & 0.355 & 0.282 & 0.217 & 0.347 & 0.251 & 0.399 & 0.350 & 0.396 & 0.084 & 0.176 \\
3,1/3 & 0.136 & 0.288 & 0.097 & 0.213 & 0.185 & 0.352 & 0.006 & 0.064 & 0.155 & 0.262 & 0.073 & 0.160 \\ 
  1,1/3 & 0.227 & 0.359 & 0.154 & 0.239 & 0.380 & 0.448 & 0.000 & 0.000 & 0.666 & 0.412 & 0.041 & 0.085 \\ 
  1/2,1/3 & 0.349 & 0.418 & 0.379 & 0.299 & 0.097 & 0.215 & 0.141 & 0.230 & 0.775 & 0.244 & 0.039 & 0.146 \\ 

 \hline
 
 \text{non-MDA} & \multicolumn{2}{c|}{m=2^{2}} & \multicolumn{2}{c|}{m=2^{4}} & \multicolumn{2}{c||}{m=2^{6}} & \multicolumn{2}{c|}{m=2^{3}} & \multicolumn{2}{c|}{m=2^{6}} & \multicolumn{2}{c|}{m=2^{9}}\\ \hline
 & 0.120 & 0.282 & 0.266 & 0.442 & 0.311 & 0.462 & 0.022 & 0.136 & 0.001 & 0.032 & 0.545 & 0.498 \\

\hline
 \end{tabu}$$
}
\end{table}

\clearpage

\begin{table}
\caption{{Estimated mixing parameter $q$ and sd values of the proposed approach with ISE metric}}
{\fontsize{6pt}{6pt}\selectfont
$$\begin{tabu}[c]{|c||cc|cc|cc||cc|cc|cc|}
 \hline
 & \multicolumn{6}{c||}{n=2^8} & \multicolumn{6}{c|}{n=2^{12}} \\ \hline
 & q & {\rm sd} & q & {\rm sd} & q & {\rm sd} & q & {\rm sd} & q & {\rm sd} & q & {\rm sd}\\ \hline\hline 
 {\rm Pareto} & \multicolumn{2}{c|}{m=2^{2}} & \multicolumn{2}{c|}{m=2^{4}} & \multicolumn{2}{c||}{m=2^{6}} & \multicolumn{2}{c|}{m=2^{3}} & \multicolumn{2}{c|}{m=2^{6}} & \multicolumn{2}{c|}{m=2^{9}} \\ \hline
1/2 & 0.021 & 0.033 & 0.138 & 0.274 & 0.174 &  0.362  & 0.016 & 0.072 & 0.027 & 0.066 & 0.259 & 0.394 \\
1 & 0.020 & 0.029 & 0.050 & 0.080 & 0.414 &  0.427 & 0.007 & 0.010 & 0.024 & 0.035 & 0.207 & 0.353 \\
3 & 0.025 & 0.029 & 0.090 & 0.138 & 0.560 &  0.401  & 0.007 & 0.010 & 0.023 & 0.035 & 0.319 & 0.402 \\
10 & 0.527 & 0.218 & 0.667 & 0.108 & 0.661 &  0.337  & 0.702 & 0.014 & 0.710 & 0.037 & 0.642 & 0.244 \\

 \hline
 {\rm T} & \multicolumn{2}{c|}{m=2^{2}} & \multicolumn{2}{c|}{m=2^{4}} & \multicolumn{2}{c||}{m=2^{6}} & \multicolumn{2}{c|}{m=2^{3}} & \multicolumn{2}{c|}{m=2^{6}} & \multicolumn{2}{c|}{m=2^{9}}\\ \hline
1/2 & 0.020 & 0.030 & 0.100 & 0.208 & 0.235 &  0.390 & 0.009 & 0.014 & 0.030 & 0.085 & 0.566 &  0.427 \\ 
1/1 & 0.020 & 0.028 & 0.052 & 0.085 & 0.460 &  0.438 & 0.007 & 0.010 & 0.020 & 0.029 & 0.244 &  0.345 \\ 
3 & 0.020 & 0.025 & 0.058 & 0.086 & 0.548 &  0.414 & 0.006 & 0.008 & 0.020 & 0.027 & 0.202 &  0.321 \\ 
10 & 0.018 & 0.020 & 0.127 & 0.230 & 0.633 &  0.372 & 0.005 & 0.006 & 0.391 & 0.372 & 0.443 &  0.371 \\ 

 \hline
 {\rm Burr} & \multicolumn{2}{c|}{m=2^{2}} & \multicolumn{2}{c|}{m=2^{4}} & \multicolumn{2}{c||}{m=2^{6}} & \multicolumn{2}{c|}{m=2^{3}} & \multicolumn{2}{c|}{m=2^{6}} & \multicolumn{2}{c|}{m=2^{9}}\\ \hline
1/2,1/2  & 0.130 & 0.299 & 0.185 & 0.341 & 0.122 &  0.280 & 0.099 & 0.247 & 0.404 & 0.395 & 0.693 &  0.362 \\ 
  1,1/2  & 0.028 & 0.078 & 0.149 & 0.283 & 0.255 &  0.400 & 0.012 & 0.047 & 0.037 & 0.109 & 0.519 &  0.427 \\ 
  3,1/2  & 0.021 & 0.029 & 0.049 & 0.076 & 0.519 &  0.451 & 0.007 & 0.010 & 0.021 & 0.029 & 0.165 &  0.270 \\ 
  1/2,1  & 0.024 & 0.051 & 0.128 & 0.247 & 0.536 &  0.448 & 0.016 & 0.062 & 0.027 & 0.065 & 0.318 &  0.397 \\ 
  1,1  & 0.021 & 0.030 & 0.051 & 0.089 & 0.480 &  0.430 & 0.007 & 0.010 & 0.023 & 0.034 & 0.290 &  0.393 \\ 
  3,1  & 0.021 & 0.025 & 0.082 & 0.140 & 0.574 &  0.402 & 0.006 & 0.008 & 0.020 & 0.028 & 0.215 &  0.321 \\ 
  1/2,3  & 0.074 & 0.069 & 0.072 & 0.122 & 0.557 &  0.433 & 0.007 & 0.010 & 0.021 & 0.031 & 0.150 &  0.256 \\ 
  1,3  & 0.025 & 0.030 & 0.099 & 0.141 & 0.592 &  0.386 & 0.006 & 0.009 & 0.020 & 0.028 & 0.211 &  0.327 \\ 
  3,3  & 0.064 & 0.140 & 0.632 & 0.213 & 0.672 &  0.344 & 0.698 & 0.122 & 0.714 & 0.052 & 0.623 &  0.264 \\ 

 \hline
 {\rm Fr{e}chet} & \multicolumn{2}{c|}{m=2^{2}} & \multicolumn{2}{c|}{m=2^{4}} & \multicolumn{2}{c||}{m=2^{6}} & \multicolumn{2}{c|}{m=2^{3}} & \multicolumn{2}{c|}{m=2^{6}} & \multicolumn{2}{c|}{m=2^{9}}\\ \hline
5 & 0.080 & 0.220 & 0.227 & 0.355 & 0.478 & 0.436 & 0.128 & 0.281 & 0.217 & 0.354 & 0.429 & 0.442 \\ 
 2 & 0.025 & 0.063 & 0.120 & 0.244 & 0.504 & 0.438 & 0.015 & 0.066 & 0.024 & 0.064 & 0.364 & 0.423 \\ 
 1 & 0.020 & 0.029 & 0.049 & 0.082 & 0.472 & 0.429 & 0.008 & 0.011 & 0.021 & 0.032 & 0.233 & 0.357 \\ 
 1/2 & 0.019 & 0.027 & 0.045 & 0.075 & 0.528 & 0.445 & 0.007 & 0.010 & 0.019 & 0.027 & 0.168 & 0.284 \\ 
 1/4 & 0.019 & 0.025 & 0.045 & 0.076 & 0.423 & 0.444 & 0.005 & 0.008 & 0.018 & 0.026 & 0.196 & 0.325 \\ 
 
 \hline
 {\rm Weibull} & \multicolumn{2}{c|}{m=2^{2}} & \multicolumn{2}{c|}{m=2^{4}} & \multicolumn{2}{c||}{m=2^{6}} & \multicolumn{2}{c|}{m=2^{3}} & \multicolumn{2}{c|}{m=2^{6}} & \multicolumn{2}{c|}{m=2^{9}}\\ \hline
1/2 & 0.452 & 0.041 & 0.423 & 0.086 & 0.769 &  0.119 & 0.419 & 0.015 & 0.714 & 0.024 & 0.975 &  0.020 \\ 
1 & 0.309 & 0.057 & 0.717 & 0.265 & 0.917 &  0.090 & 0.601 & 0.274 & 0.942 & 0.010 & 0.984 &  0.006 \\ 
3 & 0.893 & 0.023 & 0.894 & 0.036 & 0.944 &  0.071 & 0.915 & 0.006 & 0.962 & 0.010 & 0.989 &  0.003 \\ 
10 & 0.930 & 0.021 & 0.908 & 0.035 & 0.952 &  0.087 & 0.937 & 0.006 & 0.969 & 0.010 & 0.989 &  0.003 \\ 
 
 \hline
 {\rm inv. Burr} & \multicolumn{2}{c|}{m=2^{2}} & \multicolumn{2}{c|}{m=2^{4}} & \multicolumn{2}{c||}{m=2^{6}} & \multicolumn{2}{c|}{m=2^{3}} & \multicolumn{2}{c|}{m=2^{6}} & \multicolumn{2}{c|}{m=2^{9}}\\ \hline
   3,2  & 0.027 & 0.034 & 0.610 & 0.349 & 0.637 &  0.397 & 0.104 & 0.267 & 0.785 & 0.038 & 0.742 &  0.228 \\ 
  1,2  & 0.041 & 0.099 & 0.782 & 0.234 & 0.756 &  0.329 & 0.836 & 0.068 & 0.835 & 0.036 & 0.782 &  0.246 \\ 
  1/2,2  & 0.481 & 0.460 & 0.854 & 0.220 & 0.761 &  0.286 & 0.898 & 0.197 & 0.921 & 0.176 & 0.805 &  0.257 \\ 
 3,1  & 0.192 & 0.364 & 0.426 & 0.439 & 0.500 &  0.453 & 0.113 & 0.296 & 0.835 & 0.046 & 0.773 &  0.248 \\ 
  1,1  & 0.335 & 0.441 & 0.768 & 0.330 & 0.738 &  0.350 & 0.364 & 0.453 & 0.925 & 0.169 & 0.802 &  0.242 \\ 
  1/2,1  & 0.488 & 0.456 & 0.830 & 0.213 & 0.736 &  0.289 & 0.855 & 0.153 & 0.950 & 0.140 & 0.775 &  0.250 \\
3,1/3  & 0.096 & 0.285 & 0.184 & 0.380 & 0.286 &  0.435 & 0.481 & 0.460 & 0.869 & 0.228 & 0.796 &  0.279 \\ 
  1,1/3  & 0.148 & 0.316 & 0.701 & 0.356 & 0.708 &  0.371 & 0.106 & 0.269 & 0.955 & 0.092 & 0.768 &  0.234 \\ 
  1/2,1/3  & 0.349 & 0.420 & 0.802 & 0.207 & 0.666 &  0.285 & 0.826 & 0.152 & 0.980 & 0.030 & 0.707 &  0.204 \\ 

 \hline
  \text{non-MDA} & \multicolumn{2}{c|}{m=2^{2}} & \multicolumn{2}{c|}{m=2^{4}} & \multicolumn{2}{c||}{m=2^{6}} & \multicolumn{2}{c|}{m=2^{3}} & \multicolumn{2}{c|}{m=2^{6}} & \multicolumn{2}{c|}{m=2^{9}}\\ \hline
& 0.025 & 0.033 & 0.170 & 0.303 & 0.378 & 0.444 & 0.006 & 0.009 & 0.695 & 0.355 & 0.431 & 0.361 \\

\hline
 \end{tabu}$$
}
\end{table}

\clearpage

 \begin{table}
 \caption{{Estimated mixing parameter $q$ and sd values of the proposed approach with AD metric}}
{ \fontsize{6pt}{6pt} \selectfont
$$ \begin{tabu}[c]{|c||cc|cc|cc||cc|cc|cc|}
  \hline
 &  \multicolumn{6}{c||}{n=2^8} &  \multicolumn{6}{c|}{n=2^{12}}  \\  \hline
 & q & { \rm sd} & q & { \rm sd} & q & { \rm sd} & q & { \rm sd} & q & { \rm sd} & q & { \rm sd} \\  \hline \hline 
 { \rm Pareto} &  \multicolumn{2}{c|}{m=2^{2}} &  \multicolumn{2}{c|}{m=2^{4}} &  \multicolumn{2}{c||}{m=2^{6}} &  \multicolumn{2}{c|}{m=2^{3}} &  \multicolumn{2}{c|}{m=2^{6}} &  \multicolumn{2}{c|}{m=2^{9}}  \\  \hline
1/2 & 0.694 & 0.389 & 0.306 & 0.272 & 0.138 &  0.271 & 0.411 & 0.426 & 0.344 & 0.263 & 0.434 &  0.331 \\ 
1 & 0.976 & 0.027 & 0.858 & 0.280 & 0.480 &  0.371 & 0.988 & 0.031 & 0.843 & 0.301 & 0.371 &  0.258 \\ 
3 & 0.990& 0.000 & 0.972 & 0.008 & 0.592 &  0.334 & 0.990& 0.000 & 0.990& 0.000 & 0.922 &  0.139 \\ 
10 & 0.990 & 0.001 & 0.971 & 0.009 & 0.714 &  0.244 & 0.990& 0.000 & 0.990& 0.000 & 0.945 &  0.059 \\ 

  \hline
 { \rm T} &  \multicolumn{2}{c|}{m=2^{2}} &  \multicolumn{2}{c|}{m=2^{4}} &  \multicolumn{2}{c||}{m=2^{6}} &  \multicolumn{2}{c|}{m=2^{3}} &  \multicolumn{2}{c|}{m=2^{6}} &  \multicolumn{2}{c|}{m=2^{9}} \\  \hline
1/2 & 0.966 & 0.138 & 0.494 & 0.398 & 0.17 &  0.289 & 0.929 & 0.234 & 0.396 & 0.322 & 0.438 &  0.342 \\ 
1 & 0.990 & 0.004 & 0.924 & 0.186 & 0.631 &  0.401 & 0.990& 0.000 & 0.974 & 0.071 & 0.400 &  0.302 \\ 
3 & 0.990& 0.000 & 0.973 & 0.008 & 0.630 &  0.355 & 0.990& 0.000 & 0.990 & 0.000 & 0.922 &  0.137 \\ 
10 & 0.990 & 0.001 & 0.972 & 0.009 & 0.660 &  0.329 & 0.990& 0.000 & 0.990 & 0.000 & 0.948 &  0.059 \\ 

  \hline
 { \rm Burr} &  \multicolumn{2}{c|}{m=2^{2}} &  \multicolumn{2}{c|}{m=2^{4}} &  \multicolumn{2}{c||}{m=2^{6}} &  \multicolumn{2}{c|}{m=2^{3}} &  \multicolumn{2}{c|}{m=2^{6}} &  \multicolumn{2}{c|}{m=2^{9}} \\  \hline
1/2,1/2  & 0.056 & 0.109 & 0.148 & 0.220 & 0.119 &  0.197 & 0.072 & 0.242 & 0.304 & 0.405 & 0.651 &  0.326 \\ 
  1,1/2  & 0.689 & 0.386 & 0.315 & 0.250 & 0.216 &  0.292 & 0.412 & 0.433 & 0.347 & 0.258 & 0.428 &  0.349 \\ 
  3,1/2  & 0.990 & 0.002 & 0.961 & 0.068 & 0.644 &  0.397 & 0.990& 0.000 & 0.983 & 0.011 & 0.681 &  0.372 \\ 
  1/2,1  & 0.645 & 0.414 & 0.409 & 0.322 & 0.556 &  0.393 & 0.439 & 0.458 & 0.357 & 0.259 & 0.437 &  0.366 \\ 
  1,1  & 0.973 & 0.045 & 0.873 & 0.264 & 0.537 &  0.364 & 0.989 & 0.005 & 0.833 & 0.309 & 0.379 &  0.266 \\ 
  3,1  & 0.990& 0.000 & 0.973 & 0.008 & 0.628 &  0.338 & 0.990& 0.000 & 0.990& 0.000 & 0.912 &  0.168 \\ 
  1/2,3  & 0.990 & 0.003 & 0.963 & 0.067 & 0.708 &  0.355 & 0.990& 0.000 & 0.983 & 0.014 & 0.729 &  0.356 \\ 
  1,3  & 0.990& 0.000 & 0.972 & 0.008 & 0.686 &  0.322 & 0.990& 0.000 & 0.990& 0.000 & 0.936 &  0.103 \\ 
  3,3  & 0.990 & 0.001 & 0.972 & 0.009 & 0.748 &  0.236 & 0.990& 0.000 & 0.990& 0.000 & 0.948 &  0.041 \\ 

  \hline
 { \rm Fr{e}chet} &  \multicolumn{2}{c|}{m=2^{2}} &  \multicolumn{2}{c|}{m=2^{4}} &  \multicolumn{2}{c||}{m=2^{6}} &  \multicolumn{2}{c|}{m=2^{3}} &  \multicolumn{2}{c|}{m=2^{6}} &  \multicolumn{2}{c|}{m=2^{9}} \\  \hline
5 & 0.077 & 0.211 & 0.205 & 0.258 & 0.464 &  0.345 & 0.069 & 0.235 & 0.117 & 0.186 & 0.375 &  0.327 \\ 
2 & 0.783 & 0.352 & 0.425 & 0.341 & 0.510 &  0.359 & 0.625 & 0.453 & 0.349 & 0.26 & 0.339 &  0.292 \\ 
1 & 0.975 & 0.039 & 0.864 & 0.275 & 0.551 &  0.360 & 0.988 & 0.031 & 0.835 & 0.311 & 0.458 &  0.315 \\ 
  1/2  & 0.990 & 0.002 & 0.967 & 0.033 & 0.623 &  0.398 & 0.990& 0.000 & 0.988 & 0.008 & 0.848 &  0.283 \\ 
1/4 & 0.990& 0.000 & 0.966 & 0.076 & 0.428 &  0.405 & 0.990& 0.000 & 0.990& 0.000 & 0.908 &  0.184 \\ 
 
  \hline
 { \rm Weibull} &  \multicolumn{2}{c|}{m=2^{2}} &  \multicolumn{2}{c|}{m=2^{4}} &  \multicolumn{2}{c||}{m=2^{6}} &  \multicolumn{2}{c|}{m=2^{3}} &  \multicolumn{2}{c|}{m=2^{6}} &  \multicolumn{2}{c|}{m=2^{9}} \\  \hline
1/2 & 0.990& 0.000 & 0.990& 0.000 & 0.985 &  0.070 & 0.990& 0.000 & 0.990& 0.000 & 0.990 &  0.000 \\ 
1 & 0.990& 0.000 & 0.990& 0.000 & 0.990 &  0.000 & 0.990& 0.000 & 0.990& 0.000 & 0.990 &  0.000 \\ 
3 & 0.990& 0.000 & 0.990& 0.000 & 0.990 &  0.000 & 0.990& 0.000 & 0.990& 0.000 & 0.990 &  0.000 \\ 
10 & 0.990& 0.000 & 0.990& 0.000 & 0.990 &  0.003 & 0.990& 0.000 & 0.990& 0.000 & 0.990 &  0.000 \\ 
 
  \hline
 { \rm inv. Burr} &  \multicolumn{2}{c|}{m=2^{2}} &  \multicolumn{2}{c|}{m=2^{4}} &  \multicolumn{2}{c||}{m=2^{6}} &  \multicolumn{2}{c|}{m=2^{3}} &  \multicolumn{2}{c|}{m=2^{6}} &  \multicolumn{2}{c|}{m=2^{9}} \\  \hline
  3,2  & 0.987 & 0.053 & 0.961 & 0.113 & 0.664 &  0.331 & 0.990& 0.000 & 0.990 & 0.001 & 0.947 &  0.065 \\ 
1,2  & 0.988 & 0.032 & 0.956 & 0.111 & 0.765 &  0.239 & 0.990& 0.000 & 0.989 & 0.001 & 0.881 &  0.149 \\ 
  1/2,2  & 0.910 & 0.212 & 0.823 & 0.196 & 0.765 &  0.179 & 0.967 & 0.111 & 0.912 & 0.145 & 0.760 &  0.170 \\ 
  3,1  & 0.705 & 0.43 & 0.786 & 0.378 & 0.499 &  0.415 & 0.973 & 0.104 & 0.990 & 0.001 & 0.910 &  0.127 \\ 
1,1  & 0.728 & 0.388 & 0.843 & 0.256 & 0.731 &  0.274 & 0.976 & 0.068 & 0.925 & 0.125 & 0.766 &  0.177 \\ 
  1/2,1  & 0.738 & 0.343 & 0.729 & 0.194 & 0.729 &  0.157 & 0.859 & 0.252 & 0.871 & 0.176 & 0.666 &  0.117 \\  
3,1/3  & 0.010 & 0.061 & 0.112 & 0.28 & 0.252 &  0.392 & 0.909 & 0.214 & 0.813 & 0.22 & 0.807 &  0.180 \\ 
  1,1/3  & 0.822 & 0.339 & 0.813 & 0.274 & 0.69 &  0.303 & 0.990 & 0.004 & 0.921 & 0.124 & 0.681 &  0.143 \\ 
  1/2,1/3  & 0.846 & 0.304 & 0.719 & 0.163 & 0.691 &  0.110 & 0.830 & 0.28 & 0.990& 0.000 & 0.666 &  0.045 \\ 

  \hline
  { \rm non-MDA} &  \multicolumn{2}{c|}{m=2^{2}} &  \multicolumn{2}{c|}{m=2^{4}} &  \multicolumn{2}{c||}{m=2^{6}} &  \multicolumn{2}{c|}{m=2^{3}} &  \multicolumn{2}{c|}{m=2^{6}} &  \multicolumn{2}{c|}{m=2^{9}} \\  \hline
 & 0.990 & 0.000 & 0.831 & 0.335 & 0.367 & 0.415 & 0.990 & 0.000 & 0.990 & 0.000 & 0.914 & 0.183  \\
 
\hline
 \end{tabu}$$
}
\end{table}

\section{Real data study}

This section describes {three} real-world case studies. The first uses Potomac River peak stream flow (cfs) data for water years (Oct--Sep) 1895--2000 at Point Rocks, Maryland. The second case uses Danish Fire Insurance data (see also \citealt{beirlant1996practical}; \citealt{moriyama2025application}). {The last is daily log returns of the S\&P-500 index}. 

In the first case study we chose a series of $n=100$ Potomac River peaks from 1901 to 2000. By maximizing the likelihood of the annual peak flows, \cite{moriyama2025application} obtained $\widehat{\gamma}_{1^*}\fallingdotseq 0.200$, where $1^*$ denotes the maximum per year. Similarly, $\widehat{\gamma}_{5^*}\fallingdotseq 0.847$, $\widehat{\gamma}_{10^*}\fallingdotseq -0.128$, and $\widehat{\gamma}_{20^*}\fallingdotseq -0.301$, where the parameters are the maxima of five years, one decade, and two decades, respectively. The estimated probabilities of the peak flow occurrence taking more than some values by the {mixture approaches} are given in Table {10}. 

{We first note all the values of the mixing ratio parameters depend on $m$. For annual or 2 decades Peak flows all the mixture estimators were nearly nonparametric. It is shown that the parametric estimator is much more pessimistic for small $m$ and becomes comparatively optimistic as the period gets longer (Table 18 in the supplementary file, which was provided by \mbox{\citealt{moriyama2025application}}). We obtained higher probabilities by the naive mixture estimators for 5 years and {some} of the proposed estimators for 1 decade.}

\begin{table}
\caption{Estimated probabilities of peak flows of the Potomac River}{\fontsize{8pt}{8pt}\selectfont
$$\begin{tabu}[c]{|c||c||ccccc|}
 \hline
 & p ~ \text{or} ~ q & 160000 & 240000 & 320000 & 400000 & 480000 \\ \hline \hline
& \multicolumn{6}{c}{\rm annual ~ Peak ~ flows} \\ \hline

 {\rm AL} & 0.000 & 0.180 & 0.070 & 0.030 & 0.020 & 0.005 \\
  {\rm CV} & 0.000 & 0.180 & 0.070 & 0.030 & 0.020 & 0.005 \\
   {\rm ISE} & 0.000 & 0.180 & 0.070 & 0.030 & 0.020 & 0.005 \\
    {\rm AD} & 0.019 & 0.183 & 0.071 & 0.031 & 0.020 & 0.005 \\
 
 \hline
& \multicolumn{6}{c}{\rm 5 ~ years ~ Peak ~ flows} \\ \hline

 {\rm AL} & 1.000 & 0.750 & 0.436 & 0.294 & 0.217 & 0.171 \\
  {\rm CV} & 1.000 & 0.750 & 0.436 & 0.294 & 0.217 & 0.171 \\
   {\rm ISE} & 0.244 & 0.659 & 0.337 & 0.179 & 0.126 & 0.060 \\
    {\rm AD} & 0.172 & 0.650 & 0.327 & 0.167 & 0.117 & 0.050 \\

 \hline
& \multicolumn{6}{c}{\rm 1 ~ decade ~ Peak ~ flows} \\ \hline

 {\rm AL} & 0.000 & 0.863 & 0.516 & 0.263 & 0.183 & 0.049 \\
  {\rm CV} & 0.000 & 0.863 & 0.516 & 0.263 & 0.183 & 0.049 \\
   {\rm ISE} & 0.615 & 0.878 & 0.577 & 0.294 & 0.146 & 0.042 \\
    {\rm AD} & 0.731 & 0.881 & 0.588 & 0.299 & 0.138 & 0.041 \\

 \hline
& \multicolumn{6}{c}{\rm 2 ~ decades ~ Peak ~ flows} \\ \hline

 {\rm AL} & 0.000 & 0.981 & 0.766 & 0.456 & 0.332 & 0.095 \\
  {\rm CV} & 0.000 & 0.981 & 0.766 & 0.456 & 0.332 & 0.095 \\
   {\rm ISE} & 0.000 & 0.981 & 0.766 & 0.456 & 0.332 & 0.095 \\
    {\rm AD} & 0.371 & 0.982 & 0.806 & 0.505 & 0.300 & 0.078 \\

 \hline
\end{tabu}$$
}
\end{table}

\begin{table}
\caption{Estimated probabilities of fire insurance loss}{\fontsize{8pt}{8pt}\selectfont
$$\begin{tabu}[c]{|c||c||ccccc|}
 \hline
 & p ~ \text{or} ~ q & ~~x_{n,5}~~ & ~~x_{n,4}~~ & ~~x_{n,3}~~ & ~~x_{n,2}~~ & ~~x_{n,1}~~ \\ \hline \hline
 
 & \multicolumn{6}{c|}{\rm 1 ~ loss} \\ \hline

 {\rm AL} & 0.805 & 0.004 & 0.003 & 0.002 & 0.001 & 0.001 \\
  {\rm CV} & 1.000 & 0.005 & 0.003 & 0.002 & 0.002 & 0.001 \\
   {\rm ISE} & 0.990 & 0.007 & 0.004 & 0.003 & 0.002 & 0.001 \\
    {\rm AD} & 0.990 & 0.007 & 0.004 & 0.003 & 0.002 & 0.001 \\

 \hline
 & \multicolumn{6}{c|}{\rm 10 ~ loss{es}} \\ \hline
 
 {\rm AL} & 0.014 & 0.014 & 0.014 & 0.005 & 0.005 & 0.002 \\
  {\rm CV} & 0.132 & 0.015 & 0.014 & 0.005 & 0.005 & 0.003 \\
   {\rm ISE} & 0.990 & 0.027 & 0.016 & 0.010 & 0.006 & 0.004 \\
    {\rm AD} & 0.990 & 0.027 & 0.016 & 0.010 & 0.006 & 0.004 \\

 \hline
 & \multicolumn{6}{c|}{\rm 30 ~ loss{es}} \\ \hline
 
 {\rm AL} & 0.000 & 0.042 & 0.042 & 0.014 & 0.014 & 0.007 \\
  {\rm CV} & 0.000 & 0.042 & 0.042 & 0.014 & 0.014 & 0.007 \\
   {\rm ISE} & 0.990 & 0.051 & 0.030 & 0.020 & 0.012 & 0.008 \\
    {\rm AD} & 0.990 & 0.051 & 0.030 & 0.020 & 0.012 & 0.008 \\

 \hline
 & \multicolumn{6}{c|}{\rm 100 ~ loss{es}} \\ \hline
 
 {\rm AL} & 1.000 & 0.151 & 0.087 & 0.058 & 0.043 & 0.033 \\
  {\rm CV} & 1.000 & 0.151 & 0.087 & 0.058 & 0.043 & 0.033 \\
   {\rm ISE} & 0.984 & 0.165 & 0.099 & 0.062 & 0.043 & 0.033 \\
    {\rm AD} & 0.989 & 0.160 & 0.097 & 0.068 & 0.049 & 0.035 \\

 \hline
 & \multicolumn{6}{c|}{\rm 200 ~ loss{es}} \\ \hline
 
 {\rm AL} & 1.000 & 0.256 & 0.150 & 0.101 & 0.074 & 0.057 \\
  {\rm CV} & 1.000 & 0.256 & 0.150 & 0.101 & 0.074 & 0.057 \\
   {\rm ISE} & 0.990 & 0.264 & 0.158 & 0.110 & 0.083 & 0.063 \\
    {\rm AD} & 0.990 & 0.264 & 0.158 & 0.110 & 0.083 & 0.063 \\

 \hline
\end{tabu}$$
}
\end{table}

\begin{table}
\caption{{Estimated probabilities of daily log returns of the S\&P-500 index}}{\fontsize{8pt}{8pt}\selectfont
$$\begin{tabu}[c]{|c||c||ccccc|}
 \hline
 & p ~ \text{or} ~ q & ~~z_{n,5}~~ & ~~z_{n,4}~~ & ~~z_{n,3}~~ & ~~z_{n,2}~~ & ~~z_{n,1}~~ \\ \hline \hline
 
 & \multicolumn{6}{c|}{\rm 1 ~ day} \\ \hline

 {\rm AL} & 0.000 & 0.006 & 0.001 & 0.001 & 0.000 & 0.000 \\
  {\rm CV} & 0.000 & 0.006 & 0.001 & 0.001 & 0.000 & 0.000 \\
   {\rm ISE} & 0.990 & 0.025 & 0.007 & 0.005 & 0.005 & 0.005 \\
    {\rm AD} & 0.990 & 0.025 & 0.007 & 0.005 & 0.005 & 0.005 \\

 \hline
 & \multicolumn{6}{c|}{\rm 10 ~ days} \\ \hline
 
 {\rm AL} & 1.000 & 0.045 & 0.014 & 0.006 & 0.003 & 0.001 \\
  {\rm CV} & 1.000 & 0.045 & 0.014 & 0.006 & 0.003 & 0.001 \\
   {\rm ISE} & 0.990 & 0.055 & 0.024 & 0.015 & 0.012 & 0.011 \\
    {\rm AD} & 0.990 & 0.055 & 0.024 & 0.015 & 0.012 & 0.011 \\

 \hline
 & \multicolumn{6}{c|}{\rm 30 ~ days} \\ \hline
 
 {\rm AL} & 1.000 & 0.085 & 0.025 & 0.009 & 0.04 & 0.002 \\
  {\rm CV} & 1.000 & 0.085 & 0.025 & 0.009 & 0.04 & 0.002 \\
   {\rm ISE} & 0.990 & 0.094 & 0.034 & 0.019 & 0.014 & 0.012 \\
    {\rm AD} & 0.990 & 0.094 & 0.034 & 0.019 & 0.014 & 0.012 \\

 \hline
 & \multicolumn{6}{c|}{\rm 100 ~ days} \\ \hline
 
 {\rm AL} & 1.000 & 0.188 & 0.074 & 0.038 & 0.022 & 0.014 \\
  {\rm CV} & 1.000 & 0.188 & 0.074 & 0.038 & 0.022 & 0.014 \\
   {\rm ISE} & 0.990 & 0.146 & 0.050 & 0.025 & 0.017 & 0.013 \\
    {\rm AD} & 0.990 & 0.146 & 0.050 & 0.025 & 0.017 & 0.013 \\

 \hline
 & \multicolumn{6}{c|}{\rm 200 ~ days} \\ \hline
 
 {\rm AL} & 1.000 & 0.325 & 0.128 & 0.061 & 0.034 & 0.021 \\
  {\rm CV} & 1.000 & 0.325 & 0.128 & 0.061 & 0.034 & 0.021 \\
   {\rm ISE} & 0.990 & 0.178 & 0.062 & 0.030 & 0.019 & 0.015 \\
    {\rm AD} & 0.990 & 0.178 & 0.062 & 0.030 & 0.019 & 0.015 \\

 \hline
\end{tabu}$$
}
\end{table}

In the second case we chose the latest $n=2100$ losses and obtained the fitted distributions. The fitted parameters are $\widehat{\gamma}_{1^*} \fallingdotseq 0.922$, where $1^*$ denotes a loss. $\widehat{\gamma}_{10^*}\fallingdotseq 0.695$, $\widehat{\gamma}_{30^*}\fallingdotseq 0.560$, $\widehat{\gamma}_{100^*}\fallingdotseq 0.718$, and $\widehat{\gamma}_{200^*}\fallingdotseq 0.706$. The estimated probabilities greater than points $x_{n,1}:=263.25$, $x_{n,2}:= (5/6)x_{n,1}$, $x_{n,3}:= (2/3)x_{n,1}$, $x_{n,4}:= (1/2)x_{n,1}$, and $x_{n,5}:= (1/3)x_{n,1}$ are summarized in Table 11.

{$\widehat{q}$ of the proposed estimators were close to one in all losses, {that is,} the mixture estimators were nearly parametric in this case. The parametric estimator returns {the} high probabilit{ies} in cases with relatively small $m$ values (Table 19 in the supplementary file, which was provided by \mbox{\citealt{moriyama2025application}}). The mixing ratio of the naive estimators also took large values for some cases; however, the ratio were numerically unstable. For 10 loss or 30 loss the mixing ratio were close to zero rather than one. In this case study the mixing ratio of the proposed mixture estimators were numerically consistent.}

{In the last case we chose the $6390$ daily log returns from $1/4/2000$ to $6/2/2025$, and considered largest drawdowns of 1 day, 10 days, 30 days, 60 days and 90 days. Then, $\widehat{\gamma}_{1^*} \fallingdotseq -0.114$, $\widehat{\gamma}_{10^*}\fallingdotseq 0.251$, $\widehat{\gamma}_{30^*}\fallingdotseq 0.228$, $\widehat{\gamma}_{60^*}\fallingdotseq 0.227$, and $\widehat{\gamma}_{90^*}\fallingdotseq 0.246$ were obtained, where $\widehat{\gamma}\fallingdotseq 0.233$ was reported in \mbox{\cite{jansen1991frequency}} in a period (see also several estimates provided in \mbox{\citealt{dehaan2006extreme}}). The estimated points were supposed to be $z_{n,1}:=0.12765$, $z_{n,2}:= (5/6)z_{n,1}$, $z_{n,3}:= (2/3)z_{n,1}$, $z_{n,4}:= (1/2)z_{n,1}$, and $z_{n,5}:= (1/3)z_{n,1}$.

Table 12 shows all the mixing ratio values were almost one except for the cases 1 day with AL or CV, where unlike the Danish Fire Insurance loss cases this result of the naive mixture approaches may be reasonable. Even if the largest drawdowns of more than 10 days approximately follow PE, each drawdown does not necessarily do that. Anyway, this case also shows the parametrically fitting estimator is considered to outperform the naive nonparametric one in the sense of describing the maxima. 

The results of the proposed mixture estimators are little different from those of the naive mixture estimators. Those of the proposed ones tend to be larger, and the difference is large especially in the case of 200 days. However, these possibly  come from the limitation on the search space (i.e. the upper limit of the proposed approaches 0.99) and almost coincide with those of the naive mixture estimators ideally.}

\section{Discussion}\label{sec12}

This section discusses properties of the mixing ratio of mixture approaches. The performance of the naive mixture approach depends on the bandwidth, which is estimated in advance. {The simulation study in section 3 revealed nontrivial difference between ALbw and CVbw. However, we note that} there are only a few bandwidth estimation methods for distribution estimations {compared with the kernel density estimation}.

The mixing ratio of the {mixture estimators} provides a good {guide to choose} between the parametric and nonparametric approaches. {However,} numerical studies showed that the mixture estimator {does not necessarily} perform as well as the naive kernel type estimator even if the estimated mixing ratio is close to zero. In such cases, it may be more beneficial to again {calculate} the naive kernel estimates{, though it results in much higher computational cost}.

In the proposed approach, there may be a better way to set the mixing parameter $q$. Concretely, this study sets $q=h(1+h)^{-1}$; however, any function $q=q(h)$ satisfying monotonicity and $0\le q\le 1$ is a candidate of the mixing parameter. {The mixing ratio of the mixture estimators depends on $m$, which is shown in the case study in Section 4. If we consider not the pointwise estimation but the global estimation in the sense of $m$, a mixing ratio being constant or smoothly varying depending on $m$ may be better and reduce the variance of the estimates as well.} This may be an interesting issue for future work.

{The result in Section {3} shows the property of the proposed mixture estimator depends on the metric. There exist various metrics other than the ISE and the Anderson-Darling-type metric measuring the distance between two probability distributions. Each metric choose the best mixing parameter in their senses. Hence, it is considered that the metric should match the quantitative evaluation of the analyst.}

{The second case study suggests that the dataset} should be estimated parametrically (i.e., not nonparametrically), which supports the findings of previous works. {In the first case study the mixing ratio is close to zero in many cases, and so the reason and the detailed examination including survey of previous works is considered to be required. It is supposed that the cases exist in which the nonparametric approach is obviously better, and so we need to continue  application study on the semiparametric mixture estimator.}

\section{Conclusion}\label{sec13}

This study proposes a computationally efficient mixture distribution estimator and applies it to SMD estimation, which is challenging and academically interesting. In SMD estimation there exist two approaches based on extreme value theory and nonparametric smoothing. The proposed mixture estimator includes the mixing ratio parameter, and numerical properties are investigated. The simulation studies demonstrate that both mixture distribution estimators can be either parametric or nonparametric. {Comparing the mixture estimators the numerical performance differs and depends on the case; however, it is reasonable for SMD estimation}. This study succeeded in demonstrating the effectiveness of the proposed mixture approach in extreme value inference. 

\section*{Acknowledgement}
The author appreciates the referees' valuable comments that helped us to improve this manuscript.

\section*{Conflict of interest}
The author declares that there are no conflicts of interest.

\section*{Data availability}
The datasets analysed during the current study are available in the \texttt{extRemes} and \texttt{evir} packages in the R software environment. The Potomac River peak stream flow (cfs) data were originally supplied by the U.S. Geological Survey (https://nwis.waterdata.usgs.gov/nwis/peak/?site\_no=01638500).

\section*{Information}
This is an Accepted Manuscript of an article published by Taylor \& Francis in Communications in Statistics -- Theory and Methods on 09 Sep 2025, available at: https://doi.org/10.1080/03610926.2025.2551071.

\bibliography{ref} 

\begin{thebibliography}{}

\bibitem[Altman and L{\'{e}}ger(1995)Altman and
  L{\'{e}}ger]{altman1995bandwidth}
Altman, N. and L{\'{e}}ger, C. (1995).
\newblock Bandwidth selection for kernel distribution function estimation.
\newblock {\em Journal of Statistical Planning and Inference\/}, {\bf 46}(2),
  195--214.

\bibitem[Beirlant {\em et~al.}(1996)Beirlant, Teugels, and
  Vynckier]{beirlant1996practical}
Beirlant, J., Teugels, J.~L., and Vynckier, P. (1996).
\newblock {\em Practical analysis of extreme values\/}, volume~50.
\newblock Leuven University Press Leuven.

\bibitem[Bowman {\em et~al.}(1998)Bowman, Hall, and Prvan]{bowman1998bandwidth}
Bowman, A., Hall, P., and Prvan, T. (1998).
\newblock {Bandwidth selection for the smoothing of distribution functions}.
\newblock {\em Biometrika\/}, {\bf 85}(4), 799--808.

\bibitem[De~Haan and Ferreira(2006)De~Haan and Ferreira]{dehaan2006extreme}
De~Haan, L. and Ferreira, A. (2006).
\newblock {\em Extreme value theory: an introduction\/}.
\newblock Springer Science \& Business Media.

\bibitem[Faraway(1990)Faraway]{faraway1990implementing}
Faraway, J. (1990).
\newblock Implementing semiparametric density estimation.
\newblock {\em Statistics \& Probability Letters\/}, {\bf 10}(2), 141--143.

\bibitem[Gbenro(2020)Gbenro]{gbenro2020using}
Gbenro, N. (2020).
\newblock Using extreme value theory to test for outliers.
\newblock {\em SSRN Electronic Journal\/}, pages 1--17.

\bibitem[Hjort and Glad(1995)Hjort and Glad]{hjort1995nonparametric}
Hjort, N.~L. and Glad, I.~K. (1995).
\newblock Nonparametric density estimation with a parametric start.
\newblock {\em The Annals of Statistics\/}, {\bf 23}(3).

\bibitem[Hjort and Jones(1996)Hjort and Jones]{hjort1996locally}
Hjort, N.~L. and Jones, M.~C. (1996).
\newblock Locally parametric nonparametric density estimation.
\newblock {\em The Annals of Statistics\/}, {\bf 24}(4), 1619--1647.

\bibitem[Jansen and De~Vries(1991)Jansen and De~Vries]{jansen1991frequency}
Jansen, D.~W. and De~Vries, C.~G. (1991).
\newblock On the frequency of large stock returns: Putting booms and busts into
  perspective.
\newblock {\em The review of economics and statistics\/}, pages 18--24.

\bibitem[Jones(1993)Jones]{jones1993kernel}
Jones, M. (1993).
\newblock Kernel density estimation when the bandwidth is large.
\newblock {\em Australian Journal of Statistics\/}, {\bf 35}(3), 319--326.

\bibitem[Kasai {\em et~al.}(2016)Kasai, Mori, Tamura, Sekine, Tsuchida, and
  Serizawa]{kasai2016predicting}
Kasai, N., Mori, S., Tamura, K., Sekine, K., Tsuchida, T., and Serizawa, Y.
  (2016).
\newblock Predicting maximum depth of corrosion using extreme value analysis
  and bayesian inference.
\newblock {\em International Journal of Pressure Vessels and Piping\/}, {\bf
  146}.

\bibitem[Komukai and Kasahara(1994)Komukai and
  Kasahara]{komukai1994requirements}
Komukai, S. and Kasahara, K. (1994).
\newblock On the requirements for a reasonable extreme-value prediction of
  maximum pits on hot-water-supply copper tubing.
\newblock {\em Journal of Research of the National Institute of Standards and
  Technology\/}, {\bf 99}(4), 321.

\bibitem[Loader(1996)Loader]{loader1996local}
Loader, C.~R. (1996).
\newblock Local likelihood density estimation.
\newblock {\em The Annals of Statistics\/}, {\bf 24}(4).

\bibitem[Mittnik {\em et~al.}(2001)Mittnik, Rachev, and
  Samorodnitsky]{mittnik2001distribution}
Mittnik, S., Rachev, S., and Samorodnitsky, G. (2001).
\newblock The distribution of test statistics for outlier detection in
  heavy-tailed samples.
\newblock {\em Mathematical and Computer Modelling\/}, {\bf 34}(9), 1171--1183.

\bibitem[Moriyama(2025)Moriyama]{moriyama2025application}
Moriyama, T. (2025).
\newblock Application of nonparametric approach to extreme value inference in
  distribution estimation of sample maximum and its properties.
\newblock {\em Australian \& New Zealand Journal of Statistics\/}, {\bf 67}(1),
  51--76.

\bibitem[Olkin and Spiegelman(1987)Olkin and
  Spiegelman]{olkin1987semiparametric}
Olkin, I. and Spiegelman, C.~H. (1987).
\newblock A semiparametric approach to density estimation.
\newblock {\em Journal of the American Statistical Association\/}, {\bf
  82}(399), 858--865.

\bibitem[Rahman {\em et~al.}(1997)Rahman, Beaver, and Gokhale]{rahman1997note}
Rahman, M., Beaver, R.~J., and Gokhale, D.~V. (1997).
\newblock A note on estimating the combining constant in semiparametric density
  estimation.
\newblock {\em Brazilian Journal of Probability and Statistics\/}, {\bf 11}(1),
  37--50.

\bibitem[Resnick(1997)Resnick]{resnick1997discussion}
Resnick, S.~I. (1997).
\newblock Discussion of the danish data on large fire insurance losses.
\newblock {\em {ASTIN} Bulletin\/}, {\bf 27}(1), 139--151.

\bibitem[Schuster and Yakowitz(1985)Schuster and
  Yakowitz]{schuster1985parametric}
Schuster, E. and Yakowitz, S. (1985).
\newblock Parametric/nonparametric mixture density estimation with application
  to flood-frequency analysis.
\newblock {\em JAWRA Journal of the American Water Resources Association\/},
  {\bf 21}(5), 797--804.

\bibitem[Smith(1987)Smith]{smith1987estimating}
Smith, R.~L. (1987).
\newblock {Estimating Tails of Probability Distributions}.
\newblock {\em The Annals of Statistics\/}, {\bf 15}(3), 1174--1207.

\bibitem[Soleymani and Lee(2014)Soleymani and Lee]{soleymani2014bootstrap}
Soleymani, M. and Lee, S.~M. (2014).
\newblock A bootstrap procedure for local semiparametric density estimation
  amid model uncertainties.
\newblock {\em Journal of Statistical Planning and Inference\/}, {\bf 153},
  75--86.

\bibitem[Talamakrouni {\em et~al.}(2016)Talamakrouni, Keilegom, and
  Ghouch]{talamakrouni2016parametrically}
Talamakrouni, M., Keilegom, I.~V., and Ghouch, A.~E. (2016).
\newblock Parametrically guided nonparametric density and hazard estimation
  with censored data.
\newblock {\em Computational Statistics {\&} Data Analysis\/}, {\bf 93},
  308--323.

\end{thebibliography}

\section*{Appendices: Lemma 1}

\begin{lemma}
If $x^{-1}h \to 0$ for (i) (ii) or $h \to 0$ for (iii), the MSE $\mathbb{E}[\{F^m(x) -\widehat{F}^m(x)\}^2]$ is of order
\begin{align*}
\begin{dcases}
h^{4} m^{-4\gamma} + \frac{m}{n}(1- hm^{-\gamma} ) &{\rm for ~ (i), ~ (iii)}\\
h^{4} (\ln m)^{4\kappa^{-1}(\kappa-1)} + \frac{m}{n} (1- h (\ln m)^{\kappa^{-1}(\kappa-1)} ) &{\rm for ~ (ii)}.
\end{dcases}
\end{align*}
Else if $h \to \infty$ and $x^{-1}h \to \infty$ for (i) (ii), it is of order
\begin{align*}
1 + 4^{-m} m^2 \{h^{-2} x^2 + (n h)^{-1} x\}.
\end{align*}
\end{lemma}

\begin{proof}
For the case $x^{-1}h \to 0$ for (i) (ii) or $h \to 0$ for (iii) it holds that
\begin{align*}
\mathbb{E}[F(x) - \widehat{F}_h(x)] \sim -h^2 f'(x) \int z^2 w(z) {\rm d}z 
\end{align*}
and
\begin{align*}
n \mathbb{V}[\widehat{F}_h(x)] = F(x) \{1-F(x)\} - h f(x) \int z W(z) w(z) {\rm d}z + O(h^2).
\end{align*}
The order of the MSE can be proven in the same manner as \cite{moriyama2025application}. We next focus on the case $h \to \infty$ and $x^{-1}h \to \infty$ for (i) (ii).

By the asymptotic expansions we have
\begin{align*}
\mathbb{E}[F(x) - \widehat{F}_h(x)] \sim F(x) - \frac{1}{2} - h^{-1} (x - \mu)
\end{align*}
and
\begin{align*}
n \mathbb{V}[\widehat{F}_h(x)] = h^{-1} (x- \mu) \{w(0) -1\} + O(h^{-2}) 
\end{align*}

Since
\begin{align*}
\widehat{F}^m(x) =& \exp \left(m \ln (\widehat{F}(x)) \right) = \exp \left(m\{-\ln2 +\ln (1 +2(\widehat{F}(x)-1))\} \right) \\
=& \exp\left(m \left\{-\ln2 + 2\widehat{F}(x) -1 + o_P((2\widehat{F}(x) -1)) \right\}\right),
\end{align*}
combining with
\begin{align*}
2\widehat{F}(x) -1 =& 2(\widehat{F}(x) -\mathbb{E}[\widehat{F}(x)]) + (2\mathbb{E}[\widehat{F}(x)] -1) \\
=& O_{P}((nh)^{-1/2}x^{1/2}) + O(h^{-1}x) =o_P(1).
\end{align*}
we have
\begin{align*}
\widehat{F}^m(x) &\sim \{2^{-1} + O(h^{-1}x) + O_{P}((nh)^{-1/2}x^{1/2})\}^m \\
&\sim 2^{-m} + 2^{1-m}\{O_P(m(nh)^{-1/2}x^{1/2}) + O(m h^{-1} x)\}.
\end{align*}
We have completed the proof.
\end{proof}

\section*{Appendices: Proof of Theorem 1}

\begin{proof}
The convergence rate of PE $G_{\widehat{\bm{\gamma}}}$ was obtained by \cite{moriyama2025application}. For $\ln F^m(x)$ {$=O(1)$}, $\mathbb{E}[\{F^m(x) -G_{\widehat{\bm{\gamma}}}(x)\}^2]$ is of order
$$
\begin{dcases}
m^{2\gamma\rho} + m^{-2} + n^{-1} m(m^{1+2\rho}+1) ~~~ &{\rm for ~ (i), ~ (iii)}\\
C^{-2} (\ln m)^{2} + n^{-1} m(m (\ln m)^{-2} +1) ~~~ &{\rm for ~ (ii)}.
\end{dcases}
$$

Suppose $h \to \infty$ and $x^{-1}h \to \infty$ first. Then, the order of the MSE is given by
\begin{align*}
\mathbb{E}[\{F^m(x) - G_{\widehat{\bm{\gamma}}}(x)\}^2] + (1-q)^2 \{ 1 + 4^{-m} m^{2+\gamma} (m^{\gamma} h^{-2} + n^{-1} h^{-1})\},
\end{align*}
which decreases as $h$ is large. The optimal convergence rate is same as that of PE.

For $h \to \infty$ and $x^{-1}h \to 0$ the order of the MSE is 
\begin{align*}
\mathbb{E}[\{F^m(x) - G_{\widehat{\bm{\gamma}}}(x)\}^2] + (1-q)^2 \{h^{4} m^{-4\gamma} + \frac{m}{n} (1- hm^{-\gamma}) \},
\end{align*}
which cannot be faster than that of PE. That means there are no reasons for employing the `moderately large' bandwidth satisfying $h \to \infty$ but $x^{-1}h \to 0$.

In a similar manner we see for $h \to 0$ the convergence rate is the slower of $q^2 \mathbb{E}[\{F^m(x) -G_{\widehat{\bm{\gamma}}}(x)\}^2]$ and $\mathbb{E}[\{F^m(x) - (\widehat{F}_h)^m(x)\}^2]$.
The rate is same as that of NE {or} slower than NE but faster than PE. 
\end{proof}

\end{document}